\documentclass[12pt]{article}
\usepackage{amssymb,amsmath,psfig}

\newtheorem{theorem}{Theorem}[section]
\newtheorem{lemma}{Lemma}[section]

\newtheorem{remark}{Remark}[section]

\def\no{\noindent}
\def\pa{\partial}

\def\to{\rightarrow}

\def\Box{\diamond}

\def\ds{\displaystyle}

\def\Om{\Omega}  
\def\om{\omega}
\def\ga{\gamma}
\def\eps{\varepsilon}
\def\f{\varphi}
\def\vt{\vartheta}

\newcommand{\R}{\mathbb R}

\newcommand{\into}{\int\!\!\!\int_\omega}
\newcommand{\bydef}{\stackrel{\rm def}{=}}
\let\a=\alpha

\let\g=\gamma

\let\eps=\varepsilon
\let\t=\theta
\renewcommand{\l}{\lambda}
\let\o=\omega

\def\D{{\mathrm{d}}}
\def\E{{\mathrm{e}}}
\def\I{{\mathrm{i}}}

\newcommand{\rmax}{r_{\mathrm{max}}}

\let\ls=\lesssim
\let\es=\approx
\let\ol=\overline
\newcommand{\III}{|\!|\!|}
\newcommand{\bigIII}{\big|\!\big|\!\big|}

\newcommand{\BigIII}{\Big|\!\Big|\!\Big|}

\newcommand{\HOneOne}{H{}^1_1(\o)}
\newcommand{\HOneOneZ}{{\overset{\diamond}{H}{}^1_{1}(\o)}}
\newcommand{\VOneOne}{V{}^1_1(\o)}
\newcommand{\VOneOneZ}{{\overset{\circ}{V}{}^1_{1}(\o)}}
\newcommand{\HOneMOne}{H^1_{-1}(\o)}
\newcommand{\HOneMOneZ}{{\overset{\circ}{H}{}^1_{-1}(\o)}}

\begin{document}

\centerline{\Large\bf The Fourier Singular Complement Method}

\medskip
\centerline{\Large\bf for the Poisson problem. Part II: axisymmetric domains}

\vskip1truecm

\centerline{P. Ciarlet, Jr,
\footnote{
ENSTA \& CNRS UMR 2706, 32, boulevard Victor, 75739 Paris Cedex 15, France. 
% Email: {\tt ciarlet@ensta.fr}.
{\em This author was supported in part by France/Hong Kong Joint Research
Programme}.
}
\quad 
B. Jung,
\footnote{Department of Mathematics, Chemnitz University of Technology, 
D-09107 Chemnitz, Germany.
% E-mail: {\tt beate.jung@hrz.tu-chemnitz.de}.
{\em This author was supported by DGA/DSP-ENSTA 00.60.075.00.470.75.01 Research Programme}.
}
\quad 
S. Kaddouri,
\footnote{
ENSTA \& CNRS UMR 2706, 32, boulevard Victor, 75739 Paris Cedex 15, France. 
% Email: {\tt kaddouri@ensta.fr}.
}
\quad 
S. Labrunie,
\footnote{
IECN, Universit\'e Henri Poincar\'e Nancy~I \& INRIA (projet CALVI), 54506 Vand\oe uvre-l\`es-Nancy cedex, France. 
% Email: {\tt labrunie@iecn.u-nancy.fr}.
}
\quad 
J. \ Zou
\footnote{Department of Mathematics, The Chinese University of Hong Kong, 
Shatin, N.T., Hong Kong. 
% E-mail: {\tt zou@math.cuhk.edu.hk}.
{\em The work of this author was fully supported by Hong Kong RGC grants (Project no.~403403 and CUHK4048/02P)}.
}
}
   
\abstract{\footnotesize This paper is the second part of a threefold article, aimed at solving 
numerically the Poisson problem in three-dimensional prismatic or axisymmetric domains.
In the first part of this series, the Fourier Singular Complement Method was introduced and analysed, in 
prismatic domains. In this second part, the FSCM is studied in axisymmetric domains with conical 
vertices, whereas, in the third part, implementation issues, numerical tests and comparisons with 
other methods are carried out. The method is based on a Fourier expansion in the direction 
parallel to the reentrant edges of the domain, and on an improved variant of the Singular 
Complement Method in the 2D section perpendicular to those edges. 
Neither refinements near the reentrant edges or vertices of the domain, nor cut-off functions are required in 
the computations to achieve an optimal convergence order in terms of the mesh size and the number 
of Fourier modes used.}

\bigskip
%{\it Mathematics Subject Classification (1991):} % 65N30, 35L15.
{\it Date of this version} :  July 7, 2005

\section{Introduction}

The {\em Singular Complement Method} (SCM) was originally introduced by Assous {\em et al.} 
\cite{AsCS98,AsCS00}, for the 2D static or instationary Maxwell equations without charges. 
It was then extended~\cite{AsCL03,ACLS03} to the \emph{fully axisymmetric} case, i.e.~axisymmetric domains \emph{and} data, with or without charges.
The SCM has been extended in \cite{CiHe03} to the 2D Poisson problem. As noted in~\cite{CJK+04a}, further extensions to the 2D heat or wave equations, or to similar problems with piecewise constant coefficients, can be obtained easily. 
Methodologically speaking, the SCM consists in adding some singular test functions to
the usual $\mathbb{P}_1$ Lagrange FEM so that it recovers the optimal $H^1$-convergence rate, even 
in non-convex domains. In the fully axisymmetric case, one may simply add one singular test function per reentrant edge, and one per conical vertex of sufficiently large aperture.
%(such a vertex is called \emph{sharp} in the sequel).

\medbreak

There exist a couple of numerical methods in the literature for accurately solving 2D Poisson problems in non-convex domains. 
The SCM is clearly different from (anisotropic) mesh refinement techniques
\cite{Raug78,Hein96}, and can be applied efficiently to instationary 
problems (see Remark~4.1 of~\cite{CJK+04a}), since it does not need the refinements of the mesh and 
thus large time steps may be allowed. However the anisotropic mesh refinement methods 
have one advantage: they require only a partial knowledge of the most singular part of 
the solution. 

\medbreak

The numerical solution of 3D singular Poisson problems is quite different from the 2D case, and much more difficult. This is a relatively new field of research: most approaches rely on anisotropic mesh refinement, see for instance~\cite{Hein96,HeNW03} and Refs. therein. 
To our knowledge, this series of papers is the first attempt to generalize the SCM for three-dimensional singular Poisson problems. 

\medbreak

The rest of the paper is organised as follows. In the next Section, 
we define the geometry of the axisymmetric domain~$\Om$, 
and the suitable framework for the study of the Poisson problem in~$\Om$ 
using a Fourier expansion with respect to the rotational angle~$\theta$, 
namely, weighted Sobolev spaces over the meridian section $\omega$.
This suggests a framework for 
building the \emph{Fourier Singular Complement Method} (FSCM) 
for accurately solving the Poisson problem, using a Fourier expansion 
in~$\theta$, and an improved variant of the Singular
Complement Method~\cite{CiHe03} in~$\omega$.
In Section~3, we study theoretically this variant of the SCM, 
based on a regular-singular splitting of the solution $u^k$ to the
2D~problem~(\ref{axi.laplace1:Delta.k}--\ref{axi.laplace1:cl}).
The main feature of the splitting is that it is chosen 
\emph{independently of the Fourier index~$k$} as soon as
$|k|\ge2$; this independence is important, and very helpful, from the 
computational point of view. 
Section~4 presents a few results of finite element theory in the
weighted Sobolev spaces.
In Section~\ref{axi.sec:SCM}, the SCM is considered from a numerical 
point of view, to approximate $u^k$ accurately, \emph{vi\=a} 
the discretization of the splitting. In the Section~\ref{axi.sec:FSCM}, 
we build the numerical algorithms which define the FSCM, and we show 
that it has the optimal convergence of order $O(h+N^{-1})$, where $h$ 
is the 2D mesh size and $N$ is the number of Fourier modes used.

\section{Poisson problem in axisymmetric domains}
\subsection{Geometric setting and notations}
\begin{figure}[h]
\centerline{\psfig{figure=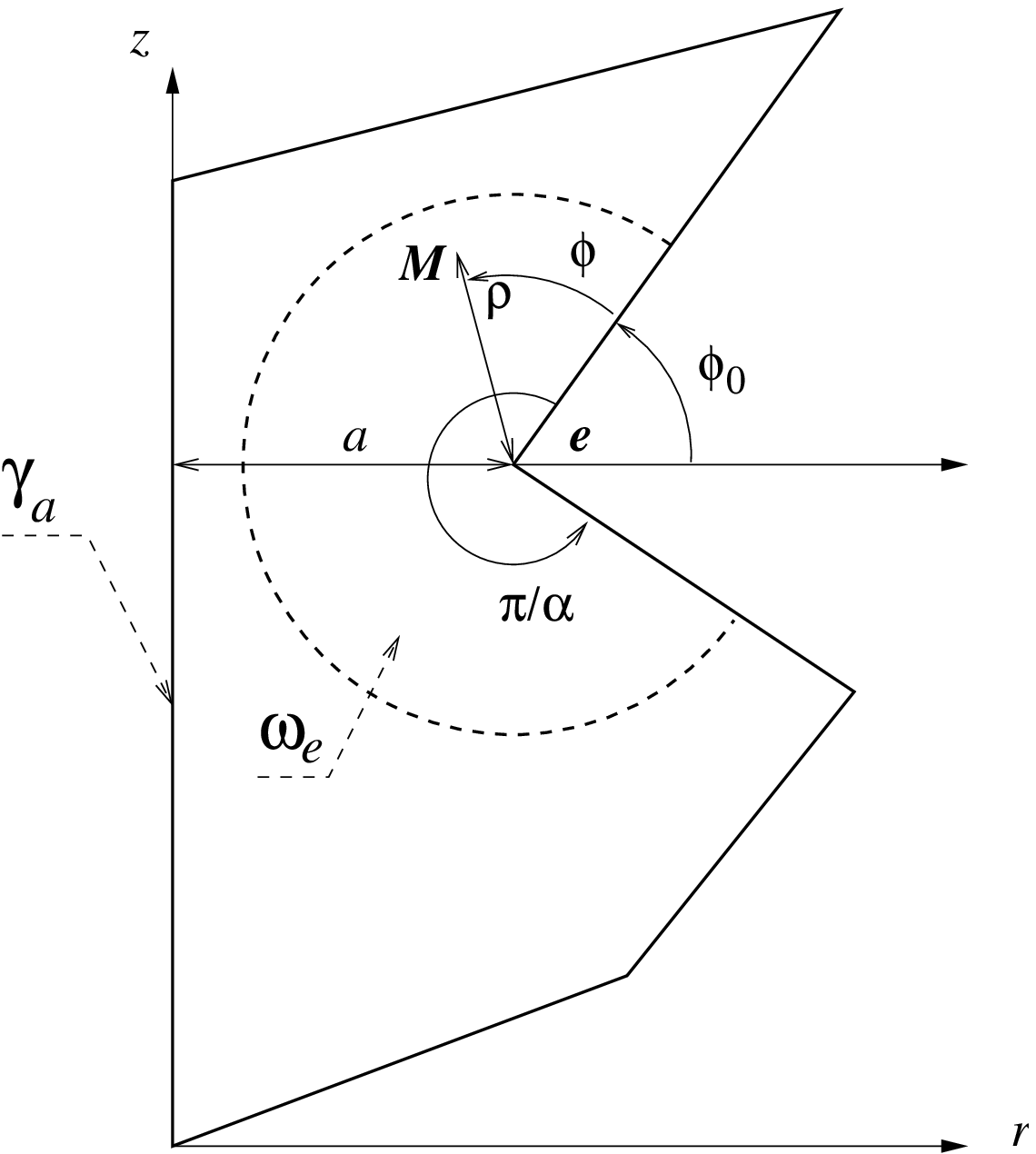,height=6.5cm}
\hspace{2cm}
\psfig{figure=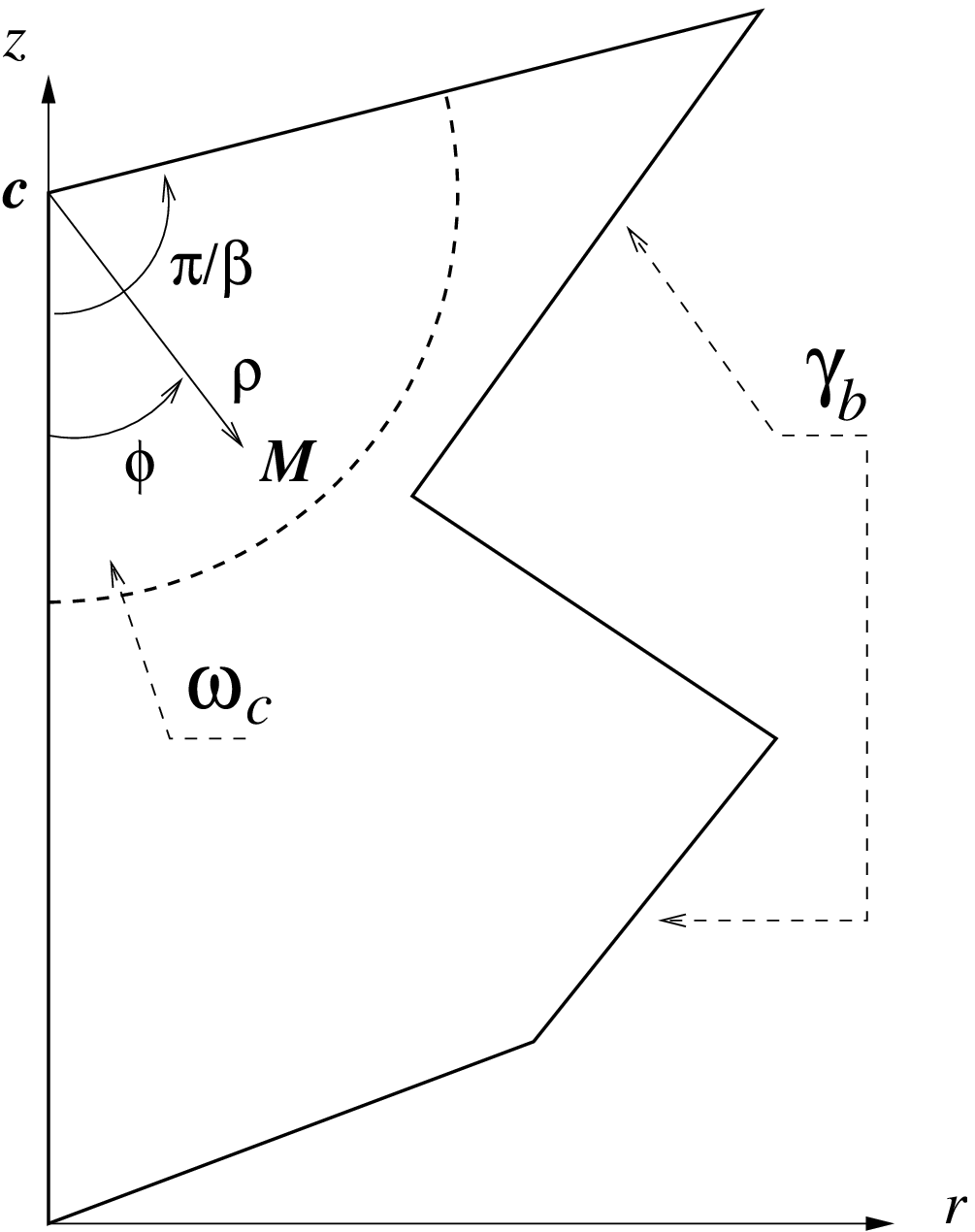,height=6.5cm} }
\caption{Notations for the geometrical singularities;
$\boldsymbol{e}$: reentrant edge; $\boldsymbol{c}$: conical vertex.}
\label{axi.fig:1}
\end{figure}

In this article, we consider an \emph{axisymmetric domain}~$\Om$, ggenerated by the 
rotation of a polygon~$\om$ around one of its sides, denoted~$\ga_a$. The 
boundary of~$\om$ is hence $\pa\om = \ga_a\cup\ga_b$, where~$\ga_b$ generates 
the boundary~$\Gamma$ of~$\Om$. Thus, $\Om$ can be described as:
\begin{equation}
 \Om =\om\times\mathbf{S}^1 \,\cup\, \ga_a. \label{axisymmetric}
\end{equation} 
The natural cylindrical coordinates will be denoted by~$(r,\theta,z)$.
The geometrical singularities that may occur on~$\Gamma$ are circular edges 
and conical vertices, which correspond to off-axis corners of~$\ga_b$ and
to its extremities. Figure~\ref{axi.fig:1} precises the various notations
associated to these singularities; a more complete description of the
geometry of~$\om$ can be found in~\cite{AsCL02,AsCL03}.

\medbreak

The problem under consideration is once more the homogeneous Dirichlet problem:
\emph{Find $u\in \overset{\circ}{H}{}^1(\Omega)$ such that}
 \begin{equation}\label{axi.Poisson}
-\Delta u = f \mbox{ in }\Om, \quad u = 0 \mbox{ on } \Gamma,
 \end{equation}
with $f\in L^2(\Om)$.
Non-homogeneous Dirichlet boundary conditions, or 
(non-)homo\-gen\-eous Neumann boundary conditions
can be handled in exactly the same manner. 

\smallbreak

As will appear in the sequel, the problem~(\ref{axi.Poisson}) will be
singular, i.e.~its solution will generically not be in~$H^2(\Om)$---as 
it would be the case in a regular or convex domain---iff
there are reentrant edges or \emph{sharp} vertices in~$\Gamma$. 
Sharp vertices are defined by the condition (see Figure~\ref{axi.fig:1}):
\begin{equation}
\nu^c < \frac12,\quad\mbox{where: } \nu^c\bydef \min\left\{ \nu>0 : P_\nu\left(\cos\frac\pi\beta \right) =0 \right\},
\end{equation}
and $P_\nu$ denotes the Legendre function. 
This is satisfied iff $\pi/\beta > \pi/\beta_\star \simeq 130^\circ48'$.
From now on, we shall assume that there is exactly \emph{one} reentrant 
edge~$\boldsymbol{e}$ (of aperture $\pi/\a$, with $1/2<\a<1$)
and \emph{one} sharp vertex~$\boldsymbol{c}$, and we shall omit the 
superscript $c$ in~$\nu^c$.

\medbreak

\paragraph{Other notations.}
We denote by $\rmax$ the supremum of the coordinate~$r$ on~$\om$, and
by $\a_0$ and $\a_1$ two fixed numbers such that 
$$1/2 < \a_0 < \a \quad\mbox{and}\quad 1/2 < \a_1 < \min(\a,\nu+1/2).$$

\medbreak

We also introduce 2D neighbourhoods $\om_e$ and $\om_c$ of 
$\boldsymbol{e}$ and~$\boldsymbol{c}$  respectively.
They stay away from all sides of~$\pa\om$ except the two ones that meet 
at the relevant corner. To them we associate cutoff functions denoted
$\eta(\rho)$, which vanish outside $\om_e$ or~$\om_c$
and depend only on the distance to the corner.

\subsection{Fourier expansions}
The functions defined on~$\Om$ will be characterised through their Fourier series in~$\theta$, the coefficients of which are functions defined on~$\om$, viz.
\begin{displaymath}
f(r,\theta,z) = \frac{1}{\sqrt{2\pi}}\, \sum_{k=-\infty}^{+\infty} f^k(r,z)\,\E^{\I k\theta},
\end{displaymath}
and the truncated Fourier expansion of $f$ at order~$N$ is:
\[
 f^{[N]}(r,\t,z) =\frac{1}{\sqrt{2\pi}}\, 
\sum_{k=-N}^N f^k(r,z) \,\E^{\I k\theta}.
\]
The regularity of the function $f$ in the scale $H^s(\Om)$ can be
characterised by that of the $(f^k)_{k\in\mathbb{Z}}$ 
in certain spaces of functions defined 
over~$\om$~\cite[\S\S{}II.1 to~II.3]{BeDM99}, namely:
$$f\in H^s(\Om),\ s\ge0 \iff \forall k\in\mathbb{Z},\ 
f^k \in H^s_{(k)}(\om),$$ 
where the $H^s_{(k)}(\om)$ are defined in turn with the help of two
different types of weighted Sobolev spaces. We shall now give these
definitions for the values of~$s$ and~$k$  chiefly needed 
in this article. The notations for the various spaces are
the same as in~\cite{BeDM99}, where the interested reader can find the 
proofs and the most general versions of the subsequent statements.

\medbreak

First, for any $\tau\in\mathbb{R}$ we consider the weighted Lebesgue space
$$L{}^2_\tau(\om) \bydef \left\{ w \mbox{ measurable on } \om : \into |w(r,z)|^2\, \,r^\tau\,\D r\,\D z < \infty \right\}.$$
This space, as well as all the spaces introduced in this article, is 
a Hermitian space of functions with \emph{complex} values. 
The scale $\left( H{}^s_\tau(\om) \right)_{s\ge0}$ is the canonical 
Sobolev scale built upon $L{}^2_\tau(\om)$, defined 
for~$s\in\mathbb{N}$ as:
$$H{}^s_\tau(\om) \bydef \left\{ w \in L{}^2_\tau(\om) : 
\pa_r^\ell \pa_z^m w \in L{}^2_\tau(\om),\ \forall\ell,m \mbox{ s.t. } 
0\le \ell+m\le s\right\},$$
and by interpolation for $s\notin\mathbb{N}$. We denote by
$\|\cdot\|_{s,\tau}$ and $|\cdot |_{s,\tau}$  the canonical norm and 
semi-norm of $H{}^s_\tau(\om)$.
% N'EST PAS UTILISÉE DANS LE TEXTE
%$\cp{\tau}$ the Poincar\'e constant in~$H{}^1_{\tau}(\om)$, i.e.
%$$\forall w\in H{}^1_{\tau}(\om)\ \mbox{ s.t. } \left.w\right|_{\g_b} 
%= 0,\quad \| w \|_{0,\tau} \le \cp{\tau} | w |_{1,\tau}.$$
%
\smallbreak

A prominent role will be played by $L{}^2_1(\om)$; its scalar product
is denoted $(\cdot|\cdot)$, without any subscript. Upon this space,
we build another, dimensionally homogeneous Sobolev scale
$\left( V{}^s_1(\om) \right)_{s\ge0}$, defined as:
$$V{}^s_1(\om) \bydef \left\{ w \in H{}^s_1(\om) : 
r^{\ell+m-s}\,\pa_r^\ell \pa_z^m w \in L{}^2_1(\om),\ 
\forall\ell,m \mbox{ s.t. } 0\le\ell+m\le \lfloor s \rfloor\right\},$$
where $\lfloor s \rfloor$~denotes the integral part of~$s$.
One can check that the general definition reduces to
\begin{displaymath}
V{}^s_1(\om) = \left\{ w \in H{}^s_1(\om) : \left. \pa_r^j w  
\right|_{\ga_a}=0, \mbox{ for } 0\le j< s-1 \right\},
\end{displaymath}
when $s$ is not an integer; while for the first values 
of~$s\in\mathbb{N}$, we have:
\begin{displaymath}
V{}^0_1(\om) = L{}^2_1(\om),\qquad 
V{}^1_1(\om) = H{}^1_1(\om) \cap L^2_{-1}(\om),\qquad
V{}^2_1(\om) = H{}^2_1(\om) \cap H^1_{-1}(\om).
\end{displaymath} 
The canonical norm of $V{}^s_1(\om)$ is denoted by 
$\III\cdot \III_{s,1}$; it is equivalent to $|\cdot|_{s,1}$ 
except for~$s\in\mathbb{N}^*$.

\medbreak

We are now ready to define the most useful spaces of Fourier coefficients.
\begin{lemma}
The spaces $H^s_{(k)}(\om)$, for $s=0,\ 1,\ 2$, are characterised as follows.
\begin{eqnarray*}
&&H^0_{(k)}(\om) = L{}^2_1(\om),\ \forall k,\quad
H^1_{(0)}(\om) = H{}^1_1(\om),\quad
H^1_{(k)}(\om) = V{}^1_1(\om),\ \forall |k|\ge 1 \,;\\
&&H^2_{(0)}(\om) = \left\{ w\in H{}^2_1(\om) : \pa_r w \in L^2_{-1}(\om) \right\},\quad
H^2_{(\pm1)}(\om) = \left\{ w\in H{}^2_1(\om) : w|_{\ga_a} = 0 \right\},\\
&&H^2_{(k)}(\om) = V{}^2_1(\om),\ \forall |k|\ge 2.
\end{eqnarray*}
\label{lem:carac:spaces}
The definition for the other values of~$s$ will be given when needed.
\end{lemma}
\begin{remark}
The scales $H{}^s_1(\om)$, $V{}^s_1(\om)$, and $H^s_{(k)}(\om)$ (for any~$k$) can be extended to negative values of the exponent~$s$, by the usual duality procedure with respect to the pivot space, which is $L{}^2_1(\om)$ in all cases.
\end{remark}

\medbreak

In order to handle the Dirichlet condition, we introduce the subspaces
$\HOneOneZ$, $\VOneOneZ$, $\HOneMOneZ$ of functions 
which vanish on~$\ga_b$. The difference in the notation is to remind that
the functions of $\VOneOne$ and $\HOneMOne\subset\VOneOne$ automatically
vanish on~$\ga_a$ in a weak sense~\cite[Prop.~4.1]{MeRa82}.
This difference is of course important when it comes to discretisation by
$\mathbb{P}_1$ finite elements.

\medbreak

Similarly to the prismatic case, we introduce the \emph{anisotropic Sobolev spaces}
\begin{eqnarray*}
h^1(\Om)&\bydef&H^1(\mathbf{S}^1,L{}^2_1(\om))
            = \{f\in L^2(\Om) : \pa_\t f\in L^2(\Om)\}~;\\
h^2(\Om)&\bydef&H^2(\mathbf{S}^1,L^2_1(\om))
            = \{f\in h^1(\Om) : \pa_{\t}^2 f\in L^2(\Om)\}~;
\end{eqnarray*}
they are identical to the $H^{0,s}(\Om)$ of~\cite[Eq.~(II.4.16)]{BeDM99},
for $s=1,\ 2$.

\medbreak

The next Lemma summarises the completeness results whose proofs
can be found in~\cite[Chapter~II]{BeDM99} or~\cite{Hein93}.
\begin{lemma}\label{axi.lem:completeness}
The following characterisations hold:
\begin{eqnarray}
f\in L^2(\Om) &\iff& \forall k\in\mathbb{Z},\ 
f^k\in L{}^2_1(\om), \mbox{ and: }
\sum_{k=-\infty}^{+\infty} \left\|f^k 
\right\|_{0,1}^2 < \infty  ~;
\label{axi.eq:fourier:L2}\\
f\in h^s(\Om) &\iff& \forall k\in\mathbb{Z},\ 
f^k\in L{}^2_1(\om), \mbox{ and: }
\sum_{k=-\infty}^{+\infty} k^{2s}\, \left\|f^k 
\right\|_{0,1}^2 < \infty ~;
\label{axi.eq:fourier:h1-2}
\end{eqnarray}
and the canonical norms of $L^2(\Om)$ and $h^s(\Om)$ are equal to the
square roots of these sums. 
Moreover, defining 
$V_{(k)} = H^1_{(k)}(\om)\cap \HOneOneZ$, viz.~$\HOneOneZ$ for $k=0$ and 
$\VOneOneZ$ for the other cases,
we have:
\begin{equation}
f\in \overset{\circ}{H}{}^1(\Omega) \iff \forall k\in\mathbb{Z},\ 
f^k\in V_{(k)} \mbox{ and }
|f|_{H^1(\Om)}^2 = \sum_{k=-\infty}^{+\infty} 
\left\|f^k \right\|_{(k)}^2 < \infty,
\label{axi.eq:fourier:H1}
\end{equation}
where the norm 
$\left\|w \right\|_{(k)}^2 = |w|_{1,1}^2 + k^2\, \|w\|_{0,-1}^2$.
\end{lemma}

\medbreak

As we did in the prismatic framework, we define the relation 
operators~$\,\ls\,$ and~$\,\es\,$ as follows. 
$a\ls b$ means $a \le C\,b$, where $C$ is a constant which depends \emph{only} 
on the geometry of the domain~$\om$, and \emph{not} on the mesh size~$h$, the 
Fourier order~$k$, or the data~$f$ of the Poisson problem. $a\es b$ denotes 
the conjunction of $a\ls b$ and $b\ls a$.

\subsection{Singular Poisson problem in 2D}
Denoting by $u^k$ and $f^k$ the Fourier coefficients of $u$ and~$f$ 
in~(\ref{axi.Poisson}), we see~\cite[\S{}II.4]{BeDM99} that for any~$k$, $u^k$ is solution to the 
following singular Poisson problem in~$\om$:\\
\emph{Find $u^k$ such that}
\begin{eqnarray}
A_k u^k \bydef -\Delta_k u^k = f^k \mbox{ in }\om, 
&&\mbox{where:}\quad \Delta_k \bydef \frac{\pa^2}{\pa r^2} + \frac1r\,\frac{\pa}{\pa r} + \frac{\pa^2}{\pa z^2} - \frac{k^2}{r^2},
\label{axi.laplace1:Delta.k}\\
u^k = 0 \mbox{ on } \g_b.&& \label{axi.laplace1:cl}
\end{eqnarray}
A special role will be played by~$\Delta_0$, whose values are the traces
in a meridian half-plane of the Laplacian of axisymmetric functions.
We remark that the operators~$\Delta_k$ have \emph{real} coefficients, hence 
the real and imaginary parts of the solution 
to~(\ref{axi.laplace1:Delta.k}--\ref{axi.laplace1:cl}) 
correspond to the real and imaginary parts of the data. So, 
\emph{in practice}, it will be sufficient to consider problems with real data 
and solutions.

\medbreak

The variational space associated 
to~(\ref{axi.laplace1:Delta.k}--\ref{axi.laplace1:cl}) is the $V_{(k)}$
defined in Lemma~\ref{axi.lem:completeness}. 
The variational formulation 
reads~\cite[\S II.4.a]{BeDM99}:
\begin{equation} \label{axi.Poissonk}
a_k\left(u^k, v\right) = \left( f^k \mid v \right), \quad \forall\, v\in 
V_{(k)},
\end{equation}
where $a_k$ is now the sesquilinear form defined by the norm 
$\|\cdot\|_{(k)}$, viz.
\begin{displaymath}\label{f4}
a_k(u,v) = \into \left[ r\,\nabla u\cdot\nabla\ol v + \frac{k^2}{r}\,u\,\ol v 
\right]\,\D\om.
\end{displaymath}
(In this text, $\nabla$~will always denote the 2D gradient in the 
$(r,z)$~plane.) 

\medbreak

Like in the prismatic case, we have the following results.
\begin{lemma}\label{axi.cor:f1} 
Let $f\in L^2(\Om)$, and $u$ be the solution to (\ref{axi.Poisson}). Then $\left( u^{[K]} \right)_K$ converges to $u$ in $H^1(\Om)$, and $\left( \Delta u^{[K]} \right)_K$ converges to $-f$ in $L^2(\Om)$.
\end{lemma}
\textit{Proof: } Similar to~\cite[Corollary~3.1]{CJK+04a}. 
\hfill$\Box$
\begin{lemma}\label{axi.cor:f2}
Let $f\in L^2(\Om)$, and $u$ be the solution to (\ref{axi.Poisson}). Then 
$\pa_\t u\in H^1(\Om)$.
\end{lemma}
\textit{Proof: } One may follow the lines of~\cite[Corollary~3.2]{CJK+04a}, 
using the \emph{a priori} estimates of~\cite[Thm~4.2]{Hein93} to check 
that $\left( \pa_r\pa_\t u,\ r^{-1}\,\pa_\t\pa_\t u,\ \pa_z\pa_\t u \right)
\in L^2(\Om)^3$.
\hfill$\Box$

\medbreak

Besides the variational space, we shall consider, for each Fourier mode~$k$:
\begin{itemize}
\item the \emph{natural} space, which is the one to which $u^k$ belongs,
i.e. the domain of the operator~$A_k$:
\begin{equation}
D(A_k) \bydef \left\{ w\in V_{(k)}=H^1_{(k)}(\om)\cap \HOneOneZ : 
A_k w\in L{}^2_1(\om) \right\} \,;
\end{equation}
\item the \emph{regularised} space, i.e. the one to which the solution $u^k$
would belong if the domain~$\Om$ were regular or convex, namely 
$H^2_{(k)}(\om)\cap \HOneOneZ$.
\end{itemize}
In~\cite[Thm~II.3.1]{BeDM99}, it is established that the regularised space no
longer depends on~$k$ as soon as $|k|\ge2$; in Theorem~\ref{axi.th:elliptic},
we will show that the same occurs for the natural space.
This suggests that the mode~2 can serve as the ``fundamental mode'' for 
the high-$|k|$ modes, just like the mode~0 does in the prismatic case.
In contradistinction to the latter, the modes 0 and~$\pm1$ have to
be treated separately, with their own singular functions.

\section{Regular-singular decompositions in the 2D domain $\om$: 
theoretical study}
We now establish the regular-singular decompositions, for the various
Fourier modes~$k$, of the solution $u^k$
to~(\ref{axi.laplace1:Delta.k}--\ref{axi.laplace1:cl}), which will be
effectively used in the numerical method. This parallels the work 
exposed in the companion paper~\cite[\S4]{CJK+04a}.

\medbreak

We shall need the following integration by parts formulae. 
\begin{theorem}\label{th.grin}
For any $u,\ v\in\HOneOne$ such that $\Delta_0 v\in L{}^2_1(\om)$, there holds:
\begin{equation}
\into \left\{ u\,\Delta_0 v + \nabla u\cdot\nabla v \right\}\,r\,\D\om = \int_{\ga_b} u\,\frac{\pa v}{\pa\nu}\, r\, \D\ga.
\label{axi.grin.1}
\end{equation}
For any~$w\in\HOneMOneZ$ such that $\Delta_0 w\in L{}^2_1(\om)$, there holds:
\begin{equation}
\Re\into -\left\{ \frac{\ol w}{r^2}\,\Delta_0 w \right\}\,r\,\D\om = \| \nabla w\|_{0,-1}^2 - 2\,\| w \|_{0,-3}^2.
\label{axi.grin.2}
\end{equation}
\end{theorem}
\textit{Proof: } Eq.~(\ref{axi.grin.1}) is the expression, in a meridian 
half-plane, of the usual Green formula applied to axisymmetric functions. 
To prove~(\ref{axi.grin.2}), we first note 
that there holds, in the sense of distributions in~$\om$:
$$\nabla w\cdot \nabla\left( \frac{\ol w}{r^2} \right) = \frac{|\nabla w|^2}{r^2} - \frac{2\,\ol w}{r^3}\,\frac{\pa w}{\pa r}.$$
But $w\in \HOneMOne$ implies $w\in L^2_{-3}(\om)$~\cite[Lemma~4.9]{AsCL02}, 
i.e. $r^{-2}\,w\in L{}^2_1(\om)$; so the above function is integrable with respect to 
the measure $r\,\D\om$, and we can apply~(\ref{axi.grin.1}) with $u=r^{-2}\,
\ol w$ and $v=w$:
\begin{eqnarray*}
I_1 &\bydef& \into -\left\{ \frac{\ol w}{r^2}\,\Delta_0 w \right\}\,r\,\D\om = 
\into \nabla w\cdot \nabla \left( \frac{\ol w}{r^2} \right)\, r\,\D\om\\
&=& \into \frac{|\nabla w|^2}{r^2}\, r\,\D\om - 2\,\into \left( \frac{1}{r}\,
\frac{\pa w}{\pa r} \right)\,  \left( \frac{\ol w}{r^2} \right)\, r\,\D\om\\
&\bydef& \| \nabla w\|_{0,-1}^2 - 2\,I_2.
\end{eqnarray*}
Now, we treat~$I_2$ by the usual integration by parts formula of order one:
\begin{eqnarray*}
\into \left[ \frac{\pa w}{\pa r}\, \frac{\ol w}{r^2} + w\,\frac{\pa}{\pa r} 
\frac{\ol w}{r^2} \right]\,\D\om = \int_{\ga} w\,\frac{\ol w}{r^2}\, \nu_r\, 
\D\om &=& 0\\
\into \left[ \frac{\pa w}{\pa r}\, \frac{\ol w}{r^2} + \frac{w}{r^2}\, 
\frac{\pa \ol w}{\pa r} + |w|^2 \left(-\frac{2}{r^3}\right) \right]\, \D\om 
&=& 0\\
2\,\Re I_2 - 2\, \| w \|_{0,-3}^2 &=& 0.
\end{eqnarray*}
Hence, $\Re I_1 =  \| \nabla w\|_{0,-1}^2 - 2\, \| w \|_{0,-3}^2$.
\hfill$\Box$

\subsection{Modes $|k|\ge2$.}
From~\cite[\S II.4]{BeDM99}, we know the following facts. 
The solution~$u^k$ to~(\ref{axi.laplace1:Delta.k}--\ref{axi.laplace1:cl}) is 
regular everywhere except in the neighbourhood of the reentrant edge, and it 
can be written as:
\begin{equation}
u^k = u^k_{_R} + \l_k\, S_k^e,\quad \mbox{with: } \left\lbrace 
\begin{matrix} 
u^k_{_R}\in H^2_{(k)}(\om)\cap\HOneOneZ = V{}^2_1(\om) \cap \VOneOneZ,\\ 
S_k^e(\rho,\phi) = \eta(\rho)\, \E^{-|k|\,\rho}\,\rho^{\a}\,\sin (\a\phi).
\hfill 
\end{matrix} \right.
\label{eq:decomp.lambda:uk}
\end{equation}
As a first consequence, we have the following
\begin{theorem}
Let $w\in D(A_k)$. Then:
\begin{itemize}
\item $w$ has a $V{}^2_1$~regularity near the axis, hence $w\in \HOneMOneZ 
\subset L^2_{-3}(\om)$, and both $\Delta_0 w$ and $r^{-2}\,w$ are in~$L{}^2_1(\om)$.
\item $w$ has an $H^{1+\a_0}$ regularity near the reentrant edge, so its 
global regularity is~$w\in V{}^{1+\a_0}_1(\om)$; and there holds: 
$\III w \III_{1+\a_0,1} \le C(k)\, \| \Delta_k w \|_{0,1}$.
\end{itemize}
\label{axi.th:elliptic}
\end{theorem}

\medbreak

In close analogy to the orthogonal decomposition of $L^2(\om)$ introduced
by Grisvard~\cite[p.~45]{Gris92}, we have:
\begin{equation}\label{axi.decomp10}
L{}^2_1(\om)=\Delta_2 [H^2_{(2)}(\om)\cap \VOneOneZ] \stackrel\perp\oplus N_2,
\end{equation}
where $N_2$ is a space of singular harmonic functions defined by 
$$N_2=\left\{p\in L^2_1(\om)\ :\ \Delta_2 p=0 \mbox{ in }\om,\ 
p=0 \mbox{ on each side of }\ga_b\right\}.$$
Here, as well as in the subsequent definitions of $N_1$ and~$N_0$,
the boundary condition on the sides of~$\ga_b$ is understood in the
suitable space, which is the trace in a meridian half-plane of the space
$\breve H(\Gamma_i)$ defined in~\cite[Definition~5.4]{AsCL02}.
Following the same line of proof as in~\cite[\S3]{AsCL03}, it is not difficult 
to establish that the dimension of $N_2$ is equal to the number of off-axis 
re-entrant corners in~$\om$, i.e.~in our case $\dim N_2=1$, and 
$N_2$=span$\{p_s^2\}$, where  $p_s^2$ can be chosen as:
\begin{equation}
p_s^2 = \mathsf{S}^e + p_{_R}^2,\quad \mbox{with: } \left\lbrace 
\begin{matrix} 
p_{_R}^2\in \VOneOneZ ,\hfill\\ 
\mathsf{S}^e(\rho,\phi) = \eta(\rho)\, \rho^{-\a}\,\sin (\a\phi).\hfill
\end{matrix} \right.
\label{eq:decomp:pSe}
\end{equation}

\medbreak

Similarly to~\cite[\S4]{CJK+04a}, we define $\f_s^2$ as the element of 
$\VOneOneZ$ which solves the Poisson problem 
\begin{equation} \label{axi.2ps}
  - \Delta_2\f_s^2 =p_s^2 \quad \mbox{in} \quad \om\,.
\end{equation}
Then by the decomposition (\ref{axi.decomp10}), 
we can split the solution $u^k$ 
to~(\ref{axi.laplace1:Delta.k}--\ref{axi.laplace1:cl}) as
\begin{equation}\label{axi.observ1}
u^k=\tilde u^k + c_k\f_s^2,
\end{equation}
where $\tilde u^k\in H^2_{(2)}(\om)\cap \VOneOneZ = 
H^2_{(k)}(\om)\cap \HOneOneZ$, and is called 
the regular part of $u^k$.  How is this decomposition related
to~(\ref{eq:decomp.lambda:uk}) ? 
Applying~(\ref{eq:decomp.lambda:uk}) to $\f_s^2$ itself gives:
$\f_s^2 = \f_{_R}^2 + \delta^2\,S_2^e$;
%\label{eq:decomp:lambda:phiSe}
observing that all the $S_k^e$ have the same principal part,
%---equal to $\rho^{\a}\,\sin (\a\phi)$---
we deduce
$u^k_{_R} = \tilde u^k + c_k\,\f_{_R}^2$, and $\l_k = c_k\,\delta^2$.
Then, using the orthogonality 
relation~(\ref{axi.decomp10}) we infer:
\begin{displaymath}
\delta^2 = \frac{ \left\| p_s^2 \right\|_{0,1}^2}{\left( -\Delta_2 S_2^e \mid p_s^2 \right) } .
\end{displaymath}
Calculating this scalar product is rather tedious but can be done 
using~(\ref{eq:decomp:pSe}) and~(\ref{axi.grin.1})---modified so as to avoid 
the singularity. We find:
\begin{equation}
\frac{\lambda_k}{c_k} = \delta^2 = \frac{1}{a\pi}\, 
\left\| p_s^2 \right\|_{0,1}^2,
\label{axi.form.beta}
\end{equation}
where $a=r(\boldsymbol{e})$ is the distance from the reentrant edge 
to the axis (see Fig.~\ref{axi.fig:1}).

\medbreak

The following lemma summarises some \emph{a priori} estimates on 
$u^k$ and~$c_k$.
\begin{lemma}\label{axi.lem:uk} 
Let $u^k$ be the solution to the Poisson problem 
(\ref{axi.laplace1:Delta.k}--\ref{axi.laplace1:cl}), then we have the 
following {\em a priori} estimates:
\begin{eqnarray}
k^2\,\|u^k\|_{0,-1}\le \rmax\,\|f^k\|_{0,1}\,, \quad 
k\,|u^k|_{1,1} &\le& \frac {\rmax}{\sqrt{2}} \|f^k\|_{0,1}\,, 
\label{axi.eq:uk1}\\
\left(k^2-2\right)\, \|u^k\|_{0,-3} \le \|f^k\|_{0,1}\,, \quad
\left(k^2-2\right)^{1/2}\, |u^k|_{1,-1} &\le& \frac {1}{\sqrt{2}} \|f^k\|_{0,1}\,,
 \label{axi.eq:uk2}\\
\|\Delta_0 u^k\|_{0,1} &\le& 2\,\|f^k\|_{0,1}\, , \label{axi.dp.uk}\\
|c_k|&\ls& k^{\alpha-1} \,\|f^k\|_{0,1}\, \label{axi.asympc}\\
\bigIII u^k \bigIII_{1+\alpha_0,1} \es
|u^k|_{1+\alpha_0,1} &\ls& k^{\alpha_0-1} \,\|f^k\|_{0,1}\,.
\label{axi.asympu}
\end{eqnarray}
\end{lemma}
\textit{Proof:} The variational formulation~(\ref{axi.Poissonk}) with 
$v=u^k$ gives:
\[
  |u^k|_{1,1}^2 + k^2 \,\|u^k\|_{0,-1}^2 \le \|f^k\|_{0,1}\,\|u^k\|_{0,1} \le 
\rmax\,\|f^k\|_{0,1}\,\|u^k\|_{0,-1}\,, 
\] 
this proves the first estimate in (\ref{axi.eq:uk1}). 
Then applying the Young inequality, we further obtain 
\[
  |u^k|_{1,1}^2 + k^2\,\|u^k\|_{0,-1}^2 \le \frac \rmax2 \left[ 
\frac{\rmax}{k^2}\, \|f^k\|_{0,1}^2 +  \frac{k^2}{\rmax}\, \|u^k\|_{0,-1}^2 
\right],
\] 
which leads to the $H{}^1_1$ semi-norm estimate in (\ref{axi.eq:uk1}). 
Similarly, multiplying~(\ref{axi.laplace1:Delta.k}) by 
$r^{-2}\,\ol u^k$ and using~(\ref{axi.grin.2}) yields:
\[
  |u^k|_{1,-1}^2 + \left(k^2-2\right) \,\|u^k\|_{0,-3}^2 \le 
\|f^k\|_{0,1}\,\|u^k\|_{0,-3} \,, 
\] 
and we obtain the two estimates in~(\ref{axi.eq:uk2}) by a similar reasoning. 
Then~(\ref{axi.dp.uk}) immediately follows from 
$\Delta_0 u^k = f^k - k^2\,r^{-2}\,u^k$.

\medbreak

The formula~(\ref{axi.form.beta}) implies: $|c_k| \es |\l_k|$; thus, the 
estimate~(\ref{axi.asympc}) is clearly equivalent to: $|\l_k| \ls k^{\a-1}$. 
This, in turn, can be obtained by following the lines of 
\cite[\S2.5.2]{Gris92} or~\cite[\S5.1]{AsCL03}. As a matter of fact, the 
latter reference shows that, away from the axis, the weights in the Sobolev 
spaces and the exact form of the modified Laplacian under consideration are 
of no importance.

\medbreak

Now, setting $f^k_{_R} \bydef -\Delta_k \left( u^k - c_k\,\f_s^2 \right)$, 
i.e. 
\begin{equation}
f^k_{_R} = -\Delta_0 u^k_{_R} + \frac{k^2}{r^2}\, u^k_{_R},
\label{axi.eq:fkR}
\end{equation}
one concludes, like in the above references, that $\left\| f^k_{_R} 
\right\|_{0,1} \ls \left\| f^k \right\|_{0,1}$. Expanding the squared norm of 
the equality~(\ref{axi.eq:fkR}) and using~(\ref{axi.grin.2}) then yields:
\begin{equation}
\left\| \Delta_0  u^k_{_R} \right\|_{0,1}^2 + 2k^2\, \left| u^k_{_R} 
\right|_{1,-1}^2 + (k^4-4k^2)\,\left\|  u^k_{_R} \right\|_{0,-3}^2 \ls \left\| 
f^k \right\|_{0,1}^2 \,.
\end{equation}
On the other hand, there holds: $u^k_{_R} \in H^2_{(0)}(\om)\cap \HOneOneZ$, 
and within this space the canonical norm of $H^2_{(0)}(\om)$ is equivalent to 
the norm~$\left\| \Delta_0  w \right\|_{0,1}$ \cite[Lemma~4.7]{AsCL03}. 
So, we have both $\left| u^k_{_R} 
\right|_{2,1} \ls \left\| f^k \right\|_{0,1}$ and  $\left| u^k_{_R} \right|_{1,1} \ls k^{-1}\, \left\| f^k \right\|_{0,1}$; and we obtain by interpolating in the scale $H{}^s_1(\om)$ that:
$ \left| u^k_{_R} \right|_{1+\a_0,1} \ls k^{\a_0-1}\, \left\| f^k \right\|_{0,1}$.
We then derive~(\ref{axi.asympu}) by adapting the proof of Lemma~4.1 of~\cite{CJK+04a}.
\hfill$\Box$

\medbreak

\begin{lemma}\label{axi.lem:bound1} 
The regular part~$\tilde u^k$ and the singularity coefficient~$c_k$ 
in~(\ref{axi.observ1}) are given as the unique solution of the coupled system:
\begin{eqnarray}
&&a_k(\tilde u^k, v) +c_k\,a_k(\f_s^2, v)
=\left( f^k \mid v \right), \quad \forall v\in \VOneOneZ\,, \label{axi.cmu3}\\
&&\left( \Vert p_s^2\Vert_{0,1}^2 + \mu\,\left[ \vert\f_s^2\vert_{1,-1}^2 + 
2\,\Vert\f_s^2\Vert_{0,-3}^2 \right] \right)\,c_k+
\mu\,\left(\tilde u^k \mid p_s^2 \right)= \left( f^k \mid p_s^2 \right)\,, 
\label{axi.cmu1}
\end{eqnarray}
where the symbol $\mu\bydef k^2-4$. And $\tilde u^k$ and $c_k$ have
the following stability estimates: 
\begin{displaymath}
%\|\nabla \tilde u^k\|_0+\sqrt{\mu}\, \|\tilde u^k\|_{0,1}
 \|\tilde u^k\|_{(k)}
 \ls k\, \|f^k\|_0\,,\qquad  %\label{axi.h1l2}\\
|c_k|\le 2\,\frac{\|f^k\|_{0,1}}{\|p_s^2\|_{0,1}}\,, \qquad
\bigIII \tilde u^k \bigIII_{2,1}\ls \|f^k\|_{0,1}\,. %\label{axi.oper16}
\end{displaymath}
\end{lemma}
We omit the details of the proof, which is very similar to that of 
Lemma~4.2 of~\cite{CJK+04a}. It makes use of the result:
$\left( r^{-2}\,{\f_s^2} \Bigm\vert  p_s^2 \right) = 
\vert\f_s^2\vert_{1,-1}^2 + 2\,\Vert\f_s^2\Vert_{0,-3}^2$,
which directly follows from~(\ref{axi.grin.2}) and~(\ref{axi.2ps}), as 
 $p_s^2$ and~$\f_s^2$ are \emph{real}.
The representation formula for the singularity coefficient is:
\begin{equation}
c_k =\frac{\left( f^k - ({\mu}/{r^2})\,A_k^{-1} f^k \Bigm\vert \,p_s^2 \right) 
}{\Vert p_s^2\Vert_{0,1}^2}\,. 
\label{axi.oper3}
\end{equation}
The scalar product $\left(  {r^{-2}}\,A_k^{-1} f^k \mid p_s^2 \right)$ 
in~(\ref{axi.oper3}) is defined thanks to Theorem~\ref{axi.th:elliptic}. 
We shall see---and this will be of practical relevance---that it can be 
written as $\left(  A_k^{-1} f^k \mid {r^{-2}}\,p_s^2 \right)$, i.e.~$p_s^2\in 
L^2_{-3}(\om)$. This is a consequence of the following lemma.
\begin{lemma}
The dual singularity $p_s^2$ is of $V{}^2_1$ regularity near the axis. It 
admits the following splitting:
\begin{equation}
p_s^2 = \tilde p^2 + p_{_P}^2,\quad \tilde p^2 \in\VOneOne,\quad p_{_P}^2 = 
\left( \frac{r}{a} \right)^2\, \rho^{-\a}\,\sin(\a\phi).
\label{axi.gamma0}
\end{equation}
Similarly, the primal singularity $\f_s^2$ can be represented as:
\begin{equation}
\f_s^2 = \tilde\f^2 + \delta^2\,\f_{_P}^2,\quad \tilde\f^2 \in V{}^2_1(\om),
\quad \f_{_P}^2 = \left( \frac{r}{a} \right)^2\, \rho^{\a}\,\sin(\a\phi).
%\label{axi.dec.phise}
\label{axi.gamma1}
\end{equation}
\end{lemma}
\emph{Proof: } Let $0 < a'' < a' < a$; we consider a cut-off function $\chi$ 
such that $\chi(r)=1$ for $r\le a''$ and $\chi(r)=0$ for $r\ge a'$, as well as the 
domain $\om' = \left\{ \boldsymbol{x}\in\om : r(\boldsymbol{x}) < a'  \right\}$. This 
domain has no off-axis reentrant corner (see Figure~\ref{axi.fig:1}), so there are no singularities 
of~$\Delta_2$, either primal or dual, in~$\om'$.

\medbreak

As we stay away from the reentrant corner, the splitting~(\ref{eq:decomp:pSe}) shows that $p_s^2\in V{}^1_1(\om')$. Thus, $\chi\,p_s^2\in V{}^1_1(\om')$ and it vanishes on~$\pa\om'$. Moreover:
$$\Delta_2\left(\chi\,p_s^2\right) = \chi\,\Delta_2 p_s^2 + \nabla\chi\cdot\nabla p_s^2 +  p_s^2\,\Delta\chi + \frac{p_s^2}{r}\,\frac{\pa\chi}{\pa r} \in L{}^2_1(\om'),$$
since the first term is identically zero, and the other three are smooth and vanish near the axis. We conclude from Theorem~\ref{axi.th:elliptic}, and the absence of primal singularities, that $\chi\,p_s^2\in V{}^2_1(\om')$, i.e. $p_s^2$ is~$V{}^2_1$ where $\chi=1$.

\medbreak

Now, using (see Figure~\ref{axi.fig:1}):
$$\left( \frac{r}{a} \right)^2 - 1 = \frac{(r-a)(r+a)}{a^2} = \frac{2}{a}\,\rho\,\cos(\phi+\phi_0) + \mbox{h.o.t.},$$
we remark
$$\mathsf{S}^{e} - p_{_P}^2 = \rho^{1-\a}\,(g_1(\phi) + \mbox{h.o.t.}) 
\in H^1(\om_e),\quad S_2^e - \f_{_P}^2 = \rho^{1+\a}\,(g_2(\phi) + 
\mbox{h.o.t.}) \in H^2(\om_e),$$
since the functions $g_{1,2}(\phi)$ as well as the higher-order terms 
(\mbox{h.o.t.}) are smooth.
Moreover, thanks to the factor $(r/a)^2$, $p_{_P}^2$ and $\f_{_P}^2$ are of $V{}^2_1$~regularity near the axis. The smoothness of these functions in the rest of the domain yields $\mathsf{S}^e - p_{_P}^2 \in \VOneOne$, $S_2^e - \f_{_P}^2\in V{}^2_1(\om)$. This proves~(\ref{axi.gamma0}) and~(\ref{axi.gamma1}).
\hfill$\Box$

\subsection{Modes $k=\pm1$.}
As we can see from Lemma~\ref{lem:carac:spaces}, the variational space is 
still $\VOneOneZ$; but the regularised space has changed. Once again, the 
only singularities are located at the reentrant edges. Hence, the 
solution~$u^k$ to~(\ref{axi.laplace1:Delta.k}--\ref{axi.laplace1:cl}), 
with~$k=\pm1$, can be split as:
\begin{equation}
u^k = u^k_{_R} + \l_k\, S_{\pm1}^e,\quad \mbox{with: } 
\left\lbrace  \begin{matrix} 
u^k_{_R}\in H^2_{(\pm1)}(\om)\cap\VOneOneZ = H{}^2_1(\om) \cap \VOneOneZ,\cr S_{\pm1}^e(\rho,\phi) = \eta(\rho)\, \E^{-\rho}\,\rho^{\a}\,\sin (\a\phi).\hfill
\end{matrix} \right.
\label{eq:decomp.lambda:uk:k=1}
\end{equation}
As a consequence of Theorem~\ref{axi.th:elliptic}, $\f_s^2\in D(A_1)$,
and the decomposition~(\ref{axi.observ1}) is still valid in this case. 
However, that singular function belongs to a space which appears too constrained 
for the modes~$\pm1$: it is even better decaying near the axis than the 
functions of~$H{}^2_1(\om) \cap \VOneOneZ$; moreover, this decay is lost in
the discretisation by $\mathbb{P}_1$ finite elements. 
So the representation formula~(\ref{axi.oper3}), though valid at the 
continuous level with $\mu=-3$, is numerically hardly stable and its use would
deteriorate the convergence rate of the SCM. 

\medbreak

So, it is better to use singular functions that are adapted for these modes.
Let $p_s^1$ be a basis of the dual singular space
$$N_1=\left\{p\in L^2_1(\om)\ :\ \Delta_1 p=0 \mbox{ in }\om,\  
p=0 \mbox{ on each side of }\ga_b\right\},$$
and $\f_s^1=A_1^{-1} p_s^1$. 
These functions were defined and studied in~\cite[\S4.1]{ACLS03},
and a numerical method was defined. We will introduce below 
(\S\ref{axi.sec:sg:petits:modes}) 
a slight modification of that method in order to improve the convergence
rate. For the moment, we recall that the function $u^k$ admits the splitting
\begin{equation}\label{axi.observ1:k=1}
u^k=\tilde u^k + c_k\f_s^1,
\end{equation}
where $\tilde u^k\in H^2_{(1)}(\om)\cap \HOneOneZ = H{}^2_{1}(\om)\cap\VOneOneZ$.
As we are in the ``usual'' SCM framework~\cite{CiHe03}, 
we have the simple representation formula
\begin{equation}
c_k =\frac{\left( f^k  \vert \,p_s^1 \right) }{\Vert p_s^1\Vert_{0,1}^2}\,,
\label{axi.oper3:k=1}
\end{equation}
and the regular part satisfies:
\begin{equation}
a_1(\tilde u^k, v) +c_k\,a_1(\f_s^1, v)
=\left( f^k \mid v \right) \quad \forall v\in \VOneOneZ\,, 
\label{axi.cmu3:k=1}
\end{equation}
From the above considerations easily follow the estimates:
\begin{eqnarray}
&&\left| u^k \right|_{1,1} \ls \left\| f^k \right\|_{0,1},\quad
\left\| u^k \right\|_{0,-1} \ls \left\| f^k \right\|_{0,1},\quad
\left\| \Delta_0 u^k \right\|_{0,1} \ls \left\| f^k \right\|_{0,1},\\
&&\left| u^k \right|_{1+\a_0,1} \ls \left\| f^k \right\|_{0,1},\quad
\bigIII \tilde u^k \bigIII_{1,1} \ls \left\| f^k \right\|_{0,1},\quad
\left| c_k \right| \ls \left\| f^k \right\|_{0,1}.
\end{eqnarray}

\subsection{Mode $k=0$.}
Now, the variational space is $V_{(0)}=\HOneOneZ$, and the regularised space is
$H^2_{(0)}(\om)\cap V_{(0)}$, with $H^2_{(0)}(\om)$ given by 
Lemma~\ref{lem:carac:spaces}.
Moreover, there is one singularity 
per reentrant edge and one per sharp vertex, see~\cite[\S{}II.4]{BeDM99}
or~\cite[\S4.4]{AsCL03}. 
The splitting of~$u^0$ with respect to regularity thus becomes:
\begin{equation}
u^0 = u^0_{_R} + \l_0^e\, S_{0}^e + \l_0^c\, S_{0}^c,\quad \mbox{with: } \left\lbrace 
\begin{matrix}
 u^0_{_R}\in H^2_{(0)}(\om)\cap\HOneOneZ,\hfill\cr S_{0}^e(\rho,\phi) = \eta(\rho)\, \rho^{\a}\,\sin (\a\phi),\hfill\\ 
S_0^c(\rho,\phi) = \eta(\rho)\, \rho^{\nu}\,P_\nu(\cos\phi).
\end{matrix}
\right.\label{eq:decomp.lambda:uk:k=0}
\end{equation}
Once more, there holds: $H^2_{(0)}(\om)\subset H^2_{(2)}(\om)$; yet, once
more, attempting to represent the singularity at the edge with the help of 
the function~$\f_s^2$ would imperil the convergence rate of the numerical
method. As a consequence, we shall use the ``natural'' singular functions
for this mode (see~\cite{AsCL03,ACLS03} for details). The dual singular space
is
$$N_0=\left\{p\in L^2_1(\om)\ :\ \Delta_0 p=0 \mbox{ in }\om,\
p=0 \mbox{ on each side of }\ga_b\right\} \,;$$
it is of dimension two, with the basis $\left\{ p_s^{0,e} \,,\, p_s^{0,c}
\right\}$ given by
\begin{equation}
\left.\begin{matrix}
p_s^{0,e} = \mathsf{S}^e + p_{_R}^{0,e},\cr
p_s^{0,c} = \mathsf{S}^c + p_{_R}^{0,c},
\end{matrix} \right\rbrace
\quad \mbox{with:}\quad \left\lbrace 
\begin{matrix} 
p_{_R}^{0,e},\ p_{_R}^{0,c}\in \HOneOneZ ,\hfill\cr 
\mathsf{S}^c(\rho,\phi) = \eta(\rho)\, 
\rho^{-\nu-1}\,P_\nu(\cos\phi),\hfill
\end{matrix} \right.
\label{eq:decomp:pSc}
\end{equation}
and $\mathsf{S}^e$ is as in~(\ref{eq:decomp:pSe}).
The \emph{primal} singular functions~$\f_s^{0,j}\in\HOneOneZ$ ($j=e,c$)
are of course defined as: $\f_s^{0,j} = A_0^{-1}\,p_s^{0,j}$.

\begin{lemma}\label{lem:singfun.c}
The primal singular function~$\f_s^{0,j}$ admits the splitting
\begin{eqnarray}
\f_s^{0,j} &=& \f_{_R}^{0,j} + \delta^{0,j}\,S_0^j,\quad\mbox{where: }
 \f_{_R}^{0,j}\in H^2_{(0)}(\om)\cap\HOneOneZ,\mbox{ and: } 
\label{eq:decomp:lambda:phiSc}\\
\delta^{0,c} &=& \ds \left\| p_s^{0,c} \right\|_{0,1}^2 \, \left\{ (1+2\,\nu)\, \int_0^{\pi/\beta} P_\nu(\cos\phi)^2\, \sin\phi\, \D\phi 
\right\}^{-1},
\label{eq:formule:delta.0.c}\\
\delta^{0,e} &=& \ds \left\| p_s^{0,e} \right\|_{0,1}^2 / (a\pi).
\label{eq:formule:delta.0.e}
\end{eqnarray}
Equivalently, $p_s^{0,j}$ and $\f_s^{0,j}$ can be represented as
\begin{eqnarray}
p_s^{0,j} = \tilde p^{0,j} + p_{_P}^{0,j},&& 
\tilde p^{0,j} \in\HOneOne,\quad \left\lbrace\begin{matrix}
p_{_P}^{0,c} = \rho^{-\nu-1}\,P_\nu(\cos\phi),\cr
p_{_P}^{0,e} = \rho^{-\a}\,\sin(\a\phi),\hfill
\end{matrix}\right.
\label{axi.gamma0.c}\\
\f_s^{0,j} = \tilde\f^{0,j} + \delta^{0,j}\,\f_{_P}^{0,j},&& 
\tilde\f^{0,j} \in H{}^2_1(\om),\quad 
\left\lbrace\begin{matrix}
\f_{_P}^{0,c} = \rho^{\nu}\,P_\nu(\cos\phi),\cr
\f_{_P}^{0,e} = \rho^{\a}\,\sin(\a\phi).\hfill
\end{matrix}\right.
\label{axi.dec.phisc}\label{axi.gamma1.c}
\end{eqnarray}
\end{lemma}
\no\textit{Proof:}
Let us examine first the conical singularity: $j=c$.
The splitting~(\ref{eq:decomp.lambda:uk:k=0}) applied to~$\f_s^{0,c}$ 
yields:
$\f_s^{0,c} = \f_{_R}^{0,c} + \delta_e^{0,c}\,S_0^e + 
\delta_c^{0,c}\,S_0^c$.
%\label{eq:decomp:phiSc:brute}
Using the definitions of~$p_s^{0,c}$ and~$\f_s^{0,c}$ we deduce
\begin{equation}
\left\| p_s^{0,c} \right\|_{0,1}^2 = \delta_e^{0,c}\,\left( -\Delta_0 S_0^e \mid p_s^{0,c} \right) + \delta_c^{0,c}\,\left( -\Delta_0 S_0^c \mid p_s^{0,c} \right).
\label{eq:calcul:CSphiSC:e}
\end{equation}
Rewriting the first scalar product in~(\ref{eq:calcul:CSphiSC:e}) as:
$$\left( -\Delta_0 S_0^e \mid p_s^{0,c} \right) = \left( -\Delta_0 S_0^e \mid p_{_R}^{0,c} \right) + \left( -\Delta_0 S_0^e \mid \mathsf{S}^c \right),$$
we notice that the second term is zero by a support argument. To evaluate the first term, we remark that $-\Delta_0 p_{_R}^{0,c} = \Delta_0\mathsf{S}^c$ vanishing near the conical point, is smooth and belongs to~$L{}^2_1(\om)$. So, both $S_0^c$ and~$p_{_R}^{0,c}$ are functions in~$\HOneOneZ$ with Laplacian~$\Delta_0$ in~$L{}^2_1(\om)$; and we use~(\ref{axi.grin.1}) twice to obtain:
$$\left( -\Delta_0 S_0^e \mid p_{_R}^{0,c} \right) = \left( S_0^e \mid -\Delta_0 p_{_R}^{0,c} \right) = \left( S_0^e \mid \Delta_0 \mathsf{S}^c \right),$$
which again vanishes by a support argument. Finally, the last scalar product in~(\ref{eq:calcul:CSphiSC:e}) can be computed e.g.~as in~\cite{NaPl94} to obtain
\begin{displaymath}
\left( -\Delta_0 S_0^c \mid p_s^{0,c} \right) = (1+2\,\nu)\, \int_0^{\pi/\beta} P_\nu(\cos\phi)^2\, \sin\phi\, \D\phi,
%\label{eq:formule:delta.0.c}
\end{displaymath}
and~(\ref{eq:decomp:lambda:phiSc}--\ref{eq:formule:delta.0.c}) are proven. This immediately yields~(\ref{axi.gamma0.c}--\ref{axi.gamma1.c}).
Similar computations are carried out for the edge singularity; 
$\delta^{0,e}$ is computed like $\delta^2$ above.
\hfill$\Box$

\medbreak

\begin{lemma}
The solution to $-\Delta_0 u^0 = f^0$ can be represented as
\begin{equation}
u^0 = \tilde u^0 + c_{0,e}\,\f_s^{0,e} + c_{0,c}\,\f_s^{0,c},\quad\mbox{where: } \tilde u^0\in H^2_{(0)}(\om)\cap\HOneOneZ.
\label{axi.observ1:k=0}
\end{equation}
The $c_{0,j}$ are given by the representation formulae:
\begin{equation}
c_{0,j} =\frac{\left( f^0 \mid p_s^{0,j} \right)}{\Vert p_s^{0,j}\Vert_{0,1}^2}\,.
\label{axi.oper3:k=0}
\end{equation}
\end{lemma}
\no\emph{Proof.}
As the space of singularities is of dimension two, it is enough to exhibit two linearly independent functions to have a basis. This is obviously the case of $\f_s^{0,e}$ and~$\f_s^{0,c}$, which proves~(\ref{axi.observ1:k=0}). 
Taking the Laplacian~$-\Delta_0$ of this equality and the scalar product by~$p_s^{0,c}$ yields, thanks to the orthogonality property:
$$\left( f^0 \mid p_s^{0,c} \right) = c_{0,e}\,\left( -\Delta_0\f_s^{0,e} \mid p_s^{0,c} \right) + c_{0,c}\,\Vert p_s^{0,c}\Vert_{0,1}^2.$$
Then, using the decomposition~(\ref{eq:decomp:lambda:phiSc}), we obtain: 
$\left( -\Delta_0\f_s^{0,e} \mid p_s^{0,c} \right) = \delta^{0,e}\, 
\left( -\Delta_0 S_0^e \mid p_s^{0,c} \right)$, which is zero as seen in the
proof of Lemma~\ref{lem:singfun.c}. Hence~(\ref{axi.oper3:k=0}) for $j=c$;
the case $j=e$ is treated similarly.
\hfill$\Box$

\medbreak

Let us state without proof the elliptic equation satisfied by $\tilde u^0$:
\begin{equation}
a_0\left( \tilde u^0,v \right) + c_{0,e}\,a_0\left( \f_s^{0,e}, v\right) + c_{0,c}\,a_0\left( \f_s^{0,c}, v \right) = \left( f^0 \mid v\right),\quad\forall v\in \HOneOneZ,
\label{axi.cmu3:k=0}
\end{equation}
and the stability estimates on the various terms in~(\ref{axi.observ1:k=0}):
\begin{eqnarray}
\left| u^0 \right|_{1,1} \ls \left\| f^0 \right\|_{0,1},&&
\left| u^0 \right|_{1+\a_1,1} \ls \left\| f^0 \right\|_{0,1},\\
\left| \tilde u^0 \right|_{1,1} \ls \left\| f^0 \right\|_{0,1},&&
\left| c_{0,j} \right| \ls \left\| f^0 \right\|_{0,1}.
\end{eqnarray}

\section{Interpolation and projection operators}
We consider a regular triangulation of the domain~$\om$, with mesh size~$h$.
The space spanned by $\mathbb{P}_1$ finite elements on this triangulation
is denoted~$V^h$; the subspace of functions which vanish on the whole 
of~$\pa\om$ is $V^h_\circ=V^h\cap\VOneOneZ$; while 
$V^h_\diamond=V^h\cap\HOneOneZ$ is the subspace of functions which 
vanish only on~$\ga_b$. 
We introduce the usual Lagrange interpolation operator $\Pi_h$ 
as well as the weighted Cl\'ement operator $P_h$. The latter---identical 
to the operator called $\widetilde\Pi_h^0$ in~\cite[\S4]{BeBD03}---is
a local projection operator onto $\mathbb{P}_1$ in the $L{}^2_1$ sense, 
which does not take into account the nodes of the triangulation which stand 
on~$\pa\om$. Hence, it maps $\VOneOneZ$ onto~$V^h_\circ$.  

\medbreak

We now prove a few results on these operators, in the framework of weighted
Sobolev spaces of fractional order. We begin by a useful density lemma.
\begin{lemma}
$V{}^2_1(\om)\cap \VOneOneZ$ is dense within $V{}^{1+\a_0}_1(\om)\cap \VOneOneZ = H{}^{1+\a_0}_1(\om)\cap \VOneOneZ$.
\label{axi.lem:densite}
\end{lemma}
\noindent\emph{Proof: } 
Let $w\in V{}^{1+\a_0}_1(\om)\cap \VOneOneZ$ and $\eps>0$. 
The construction of $\tilde w\in V{}^2_1(\om)\cap \VOneOneZ$ such that
$\| w - \tilde w \|_{1+\a_0,1} \le\eps$ is decomposed into two steps. 
The first one will not be effectively used in this article, since we do not 
consider arbitrary functions in $V{}^{1+\a_0}_1(\om)\cap \VOneOneZ$,
but only those which belong to $D(A_2)$.

\paragraph{\sc Step 1:} 
\textsc{From \(H{}^{1+\a_0}_1(\om)\cap \VOneOneZ\) 
 to~\(D(A_2)\).}\quad
From~\cite[Thm~II.4.8]{BeDM99}, we know that $A_k^{-1}$~is an isomorphism from
$H^{s-1}_{(k)}(\om)$ to $H^{s+1}_{(k)}(\om)\cap\HOneOneZ$ for $s<\a$. 
Let $C(s,k)$ be the norm of this operator. 
Setting $g=-\Delta_2 w\in H^{\a_0-1}_{(2)}(\om) = H{}^{\a_0-1}_{1}(\om)$, 
we use the density of $H{}^0_1(\om)$ within $H{}^{\a_0-1}_{1}(\om)$ to
construct $g'\in L{}^2_1(\om)$ such that $\| g-g' \|_{1-\a_0,1} < \eps/(2C(\a_0,2))$.
Then $w'\bydef A_2^{-1} g'\in D(A_2)$ and satisfies 
$\| w-w' \|_{1+\a_0,1} < \eps/2$.

\paragraph{\sc Step 2:} 
\textsc{From $D(A_2)$ to~$V{}^2_1(\om)\cap \VOneOneZ$.}\quad
There remains to find $\tilde w\in V{}^2_1(\om)\cap \VOneOneZ$ such that
$\| w'-\tilde w\|_{1+\a_0,1} < \eps/2$. Since $D(A_2) =
V{}^2_1(\om)\cap \VOneOneZ \ \oplus \ \mathrm{span}\, S_2^e$, this is 
obviously
equivalent to find $\tilde S\in V{}^2_1(\om)\cap \VOneOneZ$ such that 
$\| S_2^e - \tilde S \|_{1+\a_0,1} \le\eps'$, for arbitrary~$\eps'$. 

\medbreak

We claim that
$\tilde S = S_2^e - S_k^e$ does the job for $k$~large enough.
As a matter of fact, $S_2^e - S_k^e = \rho^{1+\a}\,\sin(\a\phi) + 
\mbox{h.o.t.} \in H^2(\om_e)\cap 
\overset{\circ}{H}{}^1(\om_e)$; thanks to the cutoff function~$\eta$, this 
gives $S_2^e - S_k^e \in V{}^2_1(\om)\cap \VOneOneZ$.
Then, it is enough to check that
\begin{displaymath}
S_k^e \to 0 \mbox{ strongly in } H^{1+\a_0}(\om_e) \mbox{ as } k\to+\infty.
\end{displaymath}
This is done like in~\cite[Lemma~5.3.3]{Gris92}, by using the Sobolev imbedding
$H^{1+\a_0}(\om_e)\subset W^{2,p}(\om_e)$ with $p=2/(2-\a_0)$; indeed,
one calculates $\| S_k^e \|_{H^{1+\a_0}(\om_e)} \ls 
\| S_k^e \|_{W^{2,p}(\om_e)} \ls k^{2(\a_0-\a)/(2-\a_0)}$.
\hfill$\Box$

\medbreak

\begin{lemma}
For any $w\in V{}^{1+\a_0}_1(\om)\cap\VOneOneZ$, there holds:
\begin{equation}
\bigIII w - P_h w \III_{1,1} \ls h^{\a_0}\, \III w \III_{1+\a_0,1},\quad
\left\| w - P_h w \right\|_{0,1}  \ls h^{1+\a_0}\, \III w \III_{1+\a_0,1}.
\label{axi.eq:interpo:Vs1}
\end{equation}
\label{axi.lem:interpo:Vs}
\end{lemma}
\noindent\emph{Proof: } 
Assume first that $w\in V{}^2_1(\om)\cap\VOneOneZ$.
From~\cite[Thm~4.2]{BeBD03}, we know
$$h^{-1}\, \|w - P_h w\|_{0,1} + \III w - P_h w \III_{1,1} \ls 
h^{s-1}\, \III w \III_{s,1} \quad\mbox{for }s=1,\ 2\,;$$
and, from~\cite[Prop.~1.e.1]{BeDM92}, that $V{}^{1+\a_0}_1(\om)$ is the 
interpolate of order~$\a_0$ between $V{}^1_1(\om)$ and~$V{}^2_1(\om)$. 
Hence the two estimates in~(\ref{axi.eq:interpo:Vs1}) by a standard 
interpolation argument. 
Then one generalises to $w\in V{}^{1+\a_0}_1(\om)\cap\VOneOneZ$
by Lemma~\ref{axi.lem:densite}.
\hfill$\Box$

\begin{lemma}
For $w\in H{}^2_1(\om)$, there holds:
\begin{equation}
\III w - \Pi_h w \III_{1-\sigma,1} \es
\left| w - \Pi_h w \right|_{1-\sigma,1} \ls h^{1+\sigma}\, |w|_{2,1}.
\label{axi.eq:interpo:Hs1}
\end{equation}
for any $\sigma\in[0,1]$.
\label{axi.lem:interpo:Hs}
\end{lemma}
\noindent\emph{Proof: } 
It stems from~\cite[Prop.~6.1]{MeRa82} or~\cite[Prop.~4.1]{BeBD03} that
$$h^{-1}\, \| w-\Pi_h w \|_{0,1} + | w-\Pi_h w |_{1,1} \ls h\, |w|_{2,1}.$$
One concludes by interpolating in the scale~$H{}^s_1(\om)$.
\hfill$\Box$

\section{Discrete formulation, SCM\label{axi.sec:SCM}}
In \S\S\ref{axi.sec:ps} to~\ref{axi.sec:approx}, the 
superscript~$2$ in $p_s^2,\ \f_s^2,\ \delta^2$, etc. 
will generally be omitted.

\subsection{Approximation of the dual singular function $p_s^2$.%
\label{axi.sec:ps}}
We start from the decomposition~(\ref{axi.gamma0}). 
 $\tilde p$ is characterised by the three conditions
\begin{displaymath} 
\tilde p\in\VOneOne,\quad \tilde p = -p_{_P} \mbox{ on } \ga,\quad
-\Delta_2 \tilde p = \Delta_2 p_{_P} \mbox{ in } \om.
%\quad \label{axi.2ps2}
\end{displaymath} 
A direct calculation shows that, denoting $\phi'=\phi+\phi_0$ 
(see Figure~\ref{axi.fig:1}):
\begin{displaymath}
\Delta_2 p_{_P} = -\frac{5r}{a^2}\,\a\,\rho^{-\a-1}\,\sin(\a\phi+\phi').
\end{displaymath}
This function is of $H^{-1}$ regularity near the reentrant edge, and smooth 
elsewhere, so it belongs to the dual of $\VOneOneZ$. However, it should be 
noticed that $\Delta_2 p_{_P}$ \emph{never} belongs locally to~$L{}^2(\om_e)$. 
This phenomenon causes the local regularity of $\tilde p$ to be weaker than in 
the prismatic case, and dramatically deteriorates the convergence rate of the 
SCM.

\medbreak

This inconvenience can be overcome by enriching the principal part with the 
next term in the expansion of~$p_s$ near the reentrant corner. To do so, we 
look for a function in~$\VOneOne$, with a $V{}^2_1$~regularity near the axis, 
vanishing on $\pa\om_e\cap\ga$, and whose Laplacian~$\Delta_2$ is locally 
``almost equal'' to~$\Delta_2 p_{_P}$.
First, we look for a local variational solution of 
$$-\Delta Q = \a\,\rho^{-\a-1}\,\sin(\a\phi+\phi')\mbox{ in }\om_e,\quad Q = 0 
\mbox{ on } \pa\om_e\cap\ga.$$ 
By separation of variables, we obtain: $Q(\rho,\phi) = \frac{1}{2}\, 
\rho^{1-\a} \cos\phi'\, \sin(\a\phi)$. As the exponent $1-\a>0$, $Q$ does 
belong to $H^1(\om_e)$.
To obtain the $V{}^2_1$~regularity near the axis, we have to multiply it 
by~$(r/a)^2$. A simple calculation shows that:
\begin{eqnarray*}
\Delta_2 \left[ \left( \frac{r}{a} \right)^2\,Q \right] &=& -
\left( \frac{r}{a} \right)^2\, \a\,\rho^{-\a-1}\,\sin(\a\phi+\phi')\\
&& \mbox{} + \frac{5\,r}{2\,a^2}\, \rho^{-\a}\, \left[ \sin(\a\phi) - \a\,
\cos\phi'\,\sin(\a\phi + \phi') \right].
\end{eqnarray*}
Hence, the new decomposition:
\begin{equation}
p_s = p_p + \widehat p, \mbox{ where: } p_p \bydef p_{_P} - \frac{5}{a}\, 
\left( \frac{r}{a} \right)^2\, Q = \left( \frac{r}{a} \right)^2\,\left[ 1 - 
\frac{5\,\rho}{2\,a}\,\cos\phi' \right]\,\rho^{-\a}\,\sin(\a\phi),
\end{equation}
enjoys the following properties. First, $\widehat p = \tilde p + (5\,r^2/a^3)
\,Q \in\VOneOne$ and it vanishes on $\pa\om_e\cap\ga$.
Then, using
$$  \left( \frac{r}{a} \right)^2 - \frac{r}{a} = \frac{r\,(r-a)}{a^2} = 
\frac{r}{a^2}\,\rho\,\cos\phi',$$
we obtain
\begin{equation}
\vt_p \bydef \Delta_2 p_p = \frac{r}{a^3}\,\rho^{-\a}\, \left[ -\frac{25}{2}\,
\sin(\a\phi) + \frac{35}{2}\,\a\,\cos\phi'\,\sin(\a\phi + \phi')  \right] 
\in L{}^2_1(\om).
\label{axi.eq:delta2pp}
\end{equation}
As $-\Delta_2\widehat p = \vt_p$, we infer by localisation that $\widehat p 
\in H^{1+\a_0}(\om_e)$. Elsewhere, the smoothness of~$\vt_p$ implies that 
of~$\widehat p$, so $\widehat p\in V{}^{1+\a_0}_1(\om)$, and is $V{}^2_1$~near 
the axis.

\medbreak

Now, we are ready to derive the FE approximation of~$p_s$. 
The variable $\widehat p$ solves the variational problem: 
\emph{Find $\widehat p\in \VOneOne$ such that}

\begin{equation} 
\widehat p = s \mbox{ on } \pa\om,\quad\mbox{and}\quad a_2\left( \widehat p, v \right) = \left( \vt_p \mid v \right)\quad \forall v\in\VOneOneZ.
\label{axi.2ps2}
\end{equation}
Similarly to the prismatic case~\cite[\S5.1]{CJK+04a}, we introduce 
\begin{itemize}
\item the boundary function $s$ which is equal to the trace of $-p_p$, hence is
zero on the two sides that meet at the reentrant corner, and smooth elsewhere;
\item the smooth extension $\tilde s\in H{}^2_1(\om)$ of~$s$ into~$\om$;
\item the variable $p^\circ = \widehat p - \tilde s$.
\end{itemize}
In the variable~$p^\circ$, the problem~(\ref{axi.2ps2}) reads: \emph{Find $p^\circ\in \VOneOneZ$ such that}
\begin{equation}
 a_2\left( p^\circ, v \right) = \left( \vt_p \mid v \right) - a_2\left( \tilde s, v \right)\quad \forall v\in\VOneOneZ \,;\label{axi.pint1}
\end{equation}
and we have $p^\circ\in  V{}^{1+\a_0}_1(\om)\cap\VOneOneZ$.
Here, too, we approximate $\widehat p$ by $\widehat p^h = \Pi_h\tilde s + p^\circ_h$, and $p_s$ by $p_s^h = p_p + \Pi_h\tilde s + p^\circ_h$, where $p^\circ_h$ solves the approximate FE problem
\begin{equation}
 a_2\left( p^\circ_h, v_h \right) = \left( \vt_p^h \mid v_h \right) - a_2\left( \Pi_h\tilde s, v_h \right)\quad \forall v_h\in V^h_\circ.\label{axi.pint3a}
\end{equation}
The notation $\left( \vt_p^h \mid v_h \right)$ stands for an approximation by a quadrature formula of the integral $\into \vt_p(r,z)\,\ol v_h(r,z)\, r\,\D r\,\D z$, with~$\vt_p(r,z)$ given by~(\ref{axi.eq:delta2pp}). As $\vt_p\in L{}^2_1(\om)$, we can suppose that the error caused by this quadrature is bounded as
\begin{equation}
\left| \left( \vt_p^h - \vt_p \mid w_h \right) \right| \le C_\mathsf{Q}^1\,h^{\mathsf{q_1}}\, \| w_h \|_{1,1},\quad \forall w_h \in V^h,\quad \mbox{for some } C_\mathsf{Q}^1>0 \mbox{ and } \mathsf{q_1} \ge 1.
\label{axi.eq:quadra.delta2pp}
\end{equation}
This can be done e.g.~by using a sixth-order Gauss--Hammer 
formula~\cite[p.~201]{Zien77},
with seven points \emph{inside} each triangle, which does not require the
unbounded value of~$\vt_p$. Of course, if~$w_h$ vanishes on~$\ga_a$, one
can replace $\| w_h \|_{1,1}$ with the stronger norm $\| w_h \|_{(2)}$ 
in~(\ref{axi.eq:quadra.delta2pp}).

\medbreak

\begin{lemma}\label{axi.lem:bound2} 
Assume $\mathsf{q_1}\ge2$; then we have:
\begin{displaymath}
\bigIII p_s-p_s^h \bigIII_{1,1} \ls h^{\alpha_0}\,,\quad
\left\| p_s-p_s^h \right\|_{0,1} \ls h^{2\alpha_0}\,. 
\end{displaymath}
\end{lemma}
\noindent\emph{Proof: }
Subtracting~(\ref{axi.pint3a}) from~(\ref{axi.pint1}) yields:
\begin{equation}
a_2\left( p^\circ - p^\circ_h, v_h \right) = \left( \vt_p -\vt_p^h \mid v_h \right) - a_2\left( \tilde s -\Pi_h\tilde s, v_h \right)\quad \forall v_h\in V^h_\circ.
\label{axi.def:varia:p0-p0h}
\end{equation}
With $v_h = p^\circ_h - w_h$, this implies:
$$\left\| p^\circ - w_h \right\|_{(2)}^2 \ge \left\| p^\circ - p^\circ_h \right\|_{(2)}^2 + 2\,a_2\left( \tilde s -\Pi_h\tilde s, p^\circ_h - w_h \right) + 2\,\left( \vt_p -\vt_p^h \mid p^\circ_h - w_h \right).$$
Now, we set $w_h = P_h p^\circ$. Using~(\ref{axi.eq:quadra.delta2pp}), we obtain
\begin{eqnarray*}
\left\| p^\circ - p^\circ_h \right\|_{(2)}^2 &\le& \left\| p^\circ - P_h p^\circ \right\|_{(2)}^2 + 2\, \left\| \tilde s -\Pi_h\tilde s \right\|_{(2)}\, \left( \left\| p^\circ_h - p^\circ \right\|_{(2)} + \left\| p^\circ - P_h p^\circ \right\|_{(2)} \right) \\ && \mbox{} + C_\mathsf{Q}^1\,h^\mathsf{q_1}\, \left( \left\| p^\circ_h - p^\circ \right\|_{(2)} + \left\| p^\circ - P_h p^\circ \right\|_{(2)} \right).
\end{eqnarray*}
With the Young inequality, the above estimate becomes:
\begin{eqnarray*}
\left\| p^\circ - p^\circ_h \right\|_{(2)}^2 &\le& \left\| p^\circ - P_h p^\circ \right\|_{(2)}^2 + 5\,\left\| \tilde s -\Pi_h\tilde s \right\|_{(2)}^2 + \frac{1}{4}\,\left\| p^\circ_h - p^\circ \right\|_{(2)}^2 + \left\| p^\circ - P_h p^\circ \right\|_{(2)}^2 \\ && \mbox{} + \frac{C_\mathsf{Q}^1}{2}\, \left[ \left( 2\,C_\mathsf{Q}^1 + 1 \right)\, h^{2\mathsf{q_1}} + \frac{\left\| p^\circ_h - p^\circ \right\|_{(2)}^2}{2\,C_\mathsf{Q}^1}  + \left\| p^\circ - P_h p^\circ \right\|_{(2)}^2 \right].
\end{eqnarray*}
Thanks to the equivalence of norms $\|\cdot\|_{(2)} \es \III\cdot\III_{1,1}$, we are left with the estimate:
\begin{displaymath}
\bigIII p^\circ - p^\circ_h  \bigIII_{1,1}^2 \ls \bigIII p^\circ - P_h p^\circ  \bigIII_{1,1}^2 + \bigIII \tilde s -\Pi_h\tilde s \bigIII_{1,1}^2 + h^{2\mathsf{q_1}}.
\end{displaymath}
By~\cite[Prop.~6.1]{MeRa82}, there holds: $ \bigIII \tilde s -\Pi_h\tilde s \bigIII_{1,1} \ls h\, |s|_{2,1}$; by Lemma~\ref{axi.lem:interpo:Vs}, we have $\bigIII p^\circ - P_h p^\circ  \bigIII_{1,1} \ls h^{\a_0}\, \bigIII p^\circ \bigIII_{1+\a_0,1}$. $s$~and~$p^\circ$ depend only on the geometry of~$\om$ 
so all their norms can be seen as constants. Hence, as soon as $\mathsf{q_1}\ge1$, there holds:
$\bigIII p^\circ - p^\circ_h  \bigIII_{1,1}^2 \ls h^{2\a_0}$.
Finally:
\begin{equation}
\bigIII p_s - p_s^h \bigIII_{1,1} = \bigIII p^\circ + \tilde s - p^\circ_h - \Pi_h\tilde s \bigIII_{1,1} \le \bigIII p^\circ - p^\circ_h  \bigIII_{1,1} + \bigIII \tilde s -\Pi_h\tilde s \bigIII_{1,1} \ls h^{\a_0}.
\label{axi.eq:norme:pse:1.1}
\end{equation}

\medbreak

The obtention of the $L^2_1$~norm estimate also follows the prismatic case closely. Here, we define~$w$ as the variational solution in~$\VOneOneZ$ to
$$a_2\left( w,v \right) = \left( p^\circ - p^\circ_h \mid v \right),\quad \forall v\in \VOneOneZ.$$
By elliptic theory~\cite[Thm~II.4.8]{BeDM99} we know $w\in H^{1+\a_0}_{(2)}(\om) = V{}^{1+\a_0}_1(\om)$ and $\III w \III_{1+\a_0,1} \ls \left\| p^\circ - p^\circ_h \right\|_{0,1}$. Its FE approximation~$w_h$ solves
\begin{equation}
a_2\left( w_h,v_h \right) = \left( p^\circ - p^\circ_h \mid v_h \right),\quad \forall v\in V^h_\circ,
\label{axi.def:varia:wh}
\end{equation}
so $\left\| w_h \right\|_{(2)} \ls \left\| p^\circ - p^\circ_h \right\|_{0,1}$; by using C\'ea's lemma and Lemma~\ref{axi.lem:interpo:Vs}, we infer:
$$\left\| w - w_h \right\|_{(2)} \ls h^{\a_0}\, \III w \III_{1+\a_0} \ls h^{\a_0}\, \left\| p^\circ - p^\circ_h \right\|_{0,1}.$$
Then, using successively~(\ref{axi.def:varia:wh}) and~(\ref{axi.def:varia:p0-p0h}), we obtain
$$\left\| p^\circ - p^\circ_h \right\|_{0,1}^2 = a_2\left( w-w_h, p^\circ - p^\circ_h\right) +
\left( \vt_p -\vt_p^h \mid w_h \right) + a_2\left( \tilde s -\Pi_h\tilde s, w-w_h \right) - a_2(\tilde s -\Pi_h\tilde s, w).$$
This is bounded by the Cauchy inequality and~(\ref{axi.eq:quadra.delta2pp}), as well as the duality argument in the scale~$V{}^s_1(\om)$:
\begin{eqnarray*}
\left\| p^\circ - p^\circ_h \right\|_{0,1}^2 &\ls& \left\| p^\circ - p^\circ_h \right\|_{(2)}\,\left\| w - w_h \right\|_{(2)} + h^{\mathsf{q_1}}\, \left\| w_h \right\|_{(2)} \\
&& \mbox{}
+ \left\| \tilde s -\Pi_h\tilde s \right\|_{(2)}\, \left\| w - w_h \right\|_{(2)} + \III \tilde s -\Pi_h\tilde s \III_{1-\a_0,1}\, \III w \III_{1+\a_0,1}\\
&\ls& \mbox{} \left\| p^\circ - p^\circ_h \right\|_{0,1}\, \left\{ h^{\a_0}\, h^{\a_0}
 + h^{\mathsf{q}_1} + h\, \left| \tilde s \right|_{2,1}\times h^{\a_0} + h^{1+\a_0}\, \left| \tilde s \right|_{2,1}\times h^0 \right\},
\end{eqnarray*}
where we have made use of~\cite[Prop.~6.1]{MeRa82} and our 
Lemma~\ref{axi.lem:interpo:Hs}. 
In order to get the $h^{2\a_0}$ estimate, we have to suppose $\mathsf{q_1} \ge 
2$. Using once more~(\ref{axi.eq:interpo:Hs1}), we obtain:
\begin{equation}
\left\| p_s - p_s^h \right\|_{0,1} \le \left\| p^\circ - p^\circ_h 
\right\|_{0,1} + \left\| \tilde s -\Pi_h\tilde s \right\|_{0,1} \ls h^{2\a_0}.
\label{axi.eq:norme:pse:0.1}
\end{equation}
\hfill$\Box$

\medbreak

We are also confronted with the task of approximating $q_s\bydef p_s/r^2$. 
The scalar product $\left( z^k \mid q_s \right)$ (see~(\ref{axi.oper3b}) 
below) is needed to compute the singularity coefficient. However, since 
$p_s^h$ is just element-wise linear, it is locally proportional to~$r$ 
in the triangles which have one or two vertices on the axis; so 
$q_s^h\bydef r^{-2}\,p_s^h\notin L{}^2_1(\om)$. This is why we cannot hope 
to control any such thing as $\left\| q_s - q_s^h \right\|_{0,1}$.

\smallbreak

Yet, thanks to the bounds (\ref{axi.eq:norme:pse:1.1}) 
and~(\ref{axi.eq:norme:pse:0.1}) for $p_s^h - p_s$, we do have the 
weak estimates:
\begin{eqnarray}
\left| \left( q_s - q_s^h \mid v \right) \right| &\ls& h^{2\a_0}\, 
\|v\|_{0,-3}, \quad \forall v \in L^2_{-3}(\om),
\label{axi.approx:tildeq:0.-3.w}\\
\mbox{resp.}\quad
\left| \left(  q_s - q_s^h \mid v \right) \right| &\ls& h^{\a_0}\, 
\|v\|_{0,-1}, \quad \forall v \in L^2_{-1}(\om).
\label{axi.approx:tildeq:0.-1.w}
\end{eqnarray}

\subsection{Approximation of the primal singular function $\f_s^2$.%
\label{axi.sec:phis}} 
We start from~(\ref{axi.gamma1}), which is sufficient for obtaining error 
estimates similar to those of the prismatic case. Using (\ref{axi.2ps}), 
we see that $\tilde\f$, satisfying $\tilde\f=-\delta\f_{_P}$ on 
$\partial\om$, solves the variational problem: 
\begin{eqnarray}
&&a_2\left( \tilde\f, v \right) = \left(p_s \mid v \right) + \delta\,
\left(\psi_{_P} \mid v \right), \quad \forall\, v\in \VOneOneZ,
\label{axi.phi1}\\
\mbox{where:}&&
\psi_{_P}\bydef\Delta_2\f_{_P} = \frac{5\,r}{a^2}\,\a\,\rho^{\a-1}\,
\sin\left[(\a-1)\,\phi - \phi_0 \right].
\end{eqnarray}
We propose the following finite element approximation of $\tilde\f$ in $V^h$:
\[\tilde\f_h= -\delta_h\,\pi_h\f_P + \f_h^0,\]
where: $\pi_h\f_P$ is a simple lifting of the boundary condition,
cf.~\cite[Eq.~(40)]{CJK+04a};
the singularity coefficient $\delta_h$ is computed using
$\delta_h =\ds \frac{1}{a\pi} \into \left(p_s^h \right)^2\, r\,\D\om$;
and $\f_h^0\in V^h_\circ$ is such that $\tilde\f_h$ is solution to the problem: 
\begin{equation}\label{axi.beta3}
a_2\left( \tilde\f_h, v_h \right) = \left(p_s^h \mid v_h \right) + \delta_h\,\left(\psi_{_P}^h \mid v_h \right),
\quad \forall\,v_h\in V^h_\circ.
\end{equation}
Like above, we assume that the quadrature formula denoted by
$\left(\psi_{_P}^h \mid v_h \right)$ satisfies:
\begin{equation}
\left| \left( \psi_{_P}^h - \psi_{_P} \mid w_h \right) \right| \le C_\mathsf{Q}^2\,h^{\mathsf{q_2}}\, \| w_h \|_{1,1},\quad \forall w_h \in V^h,\quad \mbox{for some } C_\mathsf{Q}^2>0 \mbox{ and } \mathsf{q_2} \ge 1,
\label{axi.eq:quadra.delta2phip}
\end{equation}
where one can replace $\| w_h \|_{1,1}$ with $\| w_h \|_{(2)}$ if 
$\left.w_h\right|_{\ga_a}=0$.
Then, we propose to compute the finite element approximation of $\f_s$ as:
\[ 
 \f^h_s = \tilde\f_h +\delta_h\f_{_P}\,.
\] 

\begin{lemma}\label{axi.lem:bound3} The following error estimates hold:
\[ %\begin{equation}\label{oper10}
\III\f_s-\f_s^h\III_{1,1} \ls h\,,\quad 
\|\f_s-\f_s^h\|_{(k)} \ls k\,h\,.
\] %\end{equation}
\end{lemma}
\no\emph{Proof: } We follow the lines of the proof of Lemma~5.2 of the
companion paper~\cite{CJK+04a},
taking care of the extra error generated by the quadrature. 
Subtracting~(\ref{axi.beta3}) from~(\ref{axi.phi1}), we obtain
\begin{displaymath}
a_2\left( \tilde\f - \tilde\f_h, v_h \right) = \left(p_s - p_s^h \mid v_h 
\right) + (\delta-\delta_h)\,\left(\psi_{_P} \mid v_h \right) + \delta_h\, 
\left(\psi_{_P} - \psi_{_P}^h \mid v_h \right),
\quad \forall\,v_h\in V^h_\circ.
\end{displaymath}
So, for any 
$w_h\in V^h$ satisfying $w_h-\tilde \f_h\in V^h_\circ$:
\begin{eqnarray}
\left\|\tilde\f -\tilde\f_h \right\|_{(2)}^2&\le& \left\|\tilde\f -w_h
\right\|_{(2)}^2 +2 \Bigl\{ \left\|p_s-p_s^h \right\|_{0,1} + \left|\delta-
\delta_h \right|\,\left\|\psi_{_P}\right\|_{0,1}  \nonumber\\
&&\mbox{}\hfill + \left|\delta_h\right|\,
C_\mathsf{Q}^2\,h^\mathsf{q_2}\, \left( \|\tilde\f -\tilde\f_h\|_{(2)}+ 
\|\tilde\f -w_h\|_{(2)} \right) \Bigr\}
\nonumber\\
&\le& 2\,\left\|\tilde\f -w_h\right\|_{(2)}^2 + 
\frac 12 \left\|\tilde\f -\tilde\f_h\right\|_{(2)}^2 \nonumber\\
&& \mbox{}+C\, \left( \left\|p_s-p_s^h\right\|_{0,1}^2 + \left|\delta -
\delta_h \right|^2\,\left\|\psi_{_P}\right\|_{0,1}^2 + \left|\delta_h
\right|^2\, h^{2\mathsf{q_2}} \right).  
\label{axi.phi5}
\end{eqnarray}
But $\left\|\psi_{_P}\right\|_{0,1}$ is a constant of the domain, and the 
error on the singularity coefficient is bounded as
\begin{equation}\label{axi.gamma2}
 |\delta-\delta_h| =\frac 1{a\pi} \left|
 \left\|p_s\right\|_{0,1}^2 -\left\|p_s^h\right\|_{0,1}^2\right| \ls 
\left\|p_s-p_s^h\right\|_{0,1}\ls h^{2\alpha_0}\,, 
\end{equation}
hence $\left|\delta_h\right|\es1$. With Lemma~\ref{axi.lem:bound2}, 
(\ref{axi.phi5}) becomes
\begin{equation}
\left\|\tilde\f -\tilde\f_h \right\|_{(2)}^2 \ls \left\|\tilde\f - 
w_h\right\|_{(2)}^2 + h^{4\a_0} + h^{2\mathsf{q_2}}.\label{axi.phi6}
\end{equation}
To obtain an $h^1$ estimate, it is thus sufficient to assume 
$\mathsf{q_2}\ge1$. We then derive from~(\ref{axi.phi6}) that, 
with $w_h=\delta_h\Pi_h\tilde\f/\delta$
\begin{equation}
\left\|\tilde\f -\tilde\f_h\right\|_{(2)}^2 \ls h^{4\a_0}+ |\delta|^{-2} 
\Big\{ |\delta-\delta_h|^2\,|\tilde\f|_1^2 + 
|\delta_h|^2\,\left\|\tilde\f-\Pi_h\tilde\f\right\|_{(2)}^2 \Big\}.
\label{axi.eq:phi1}
\end{equation}
As $\tilde\f\in H{}^2_1(\om)$, we have from \cite[Prop.~6.1]{MeRa82}:
$\left\|\tilde\f-\Pi_h\tilde\f\right\|_{(2)} \es \bigIII \tilde\f-\Pi_h
\tilde\f \bigIII_{1,1} \ls h\, \left| \tilde\f \right|_{2,1}$, 
which with~(\ref{axi.gamma2}) gives: $ \bigIII \tilde\f -\tilde\f_h 
\bigIII_{1,1} \es \left\|\tilde\f -\tilde\f_h\right\|_{(2)}\ls h$, and 
finally:
\begin{eqnarray*}
\bigIII \f_s-\f_s^h \bigIII_{1,1} &\le& \bigIII \tilde\f-\tilde\f_h 
\bigIII_{1,1}+ 
\left|\delta-\delta_h\right| \,\bigIII \f_{_P} \bigIII_{1,1} \ls h.
\end{eqnarray*}
Finally, the estimate on $\left\|\f_s-\f_s^h\right\|_{(k)}$  follows from 
\[
\left\|\f_s-\f_s^h\right\|_{(k)}^2=\left\|\f_s-\f_s^h\right\|_{(2)}^2 + 
\mu\,\left\|\f_s-\f_s^h\right\|_{0,-1}^2  \le\left(1 + \frac\mu4 \right)\,
\left\|\f_s-\f_s^h\right\|_{(2)}^2 \ls k^2\,h^2 .
\]
\hfill $\Box$

\subsection{Approximation of $\tilde u^k$ and $c_k$ in decomposition 
(\ref{axi.observ1}), for $|k|\ge2$.\label{axi.sec:approx}}  
Noting that $\tilde u^k$ and $c_k$ solve the coupled system 
(\ref{axi.cmu3}--\ref{axi.cmu1}),  it seems natural 
to formulate their finite element approximations as follows:\\
\emph{Find  $\tilde u^k_h\in V^h_\circ$ and $c^h_k\in \R^1$ such that:}
\begin{eqnarray}
&&a_k(\tilde u^k_h, v_h) +c_k^h\,a_k(\f_s^h, v_h)
=\left( f^k \mid v_h \right) \quad \forall v_h\in V^h_\circ\,, \label{axi.cmu6}\\
&&\left( \Vert p_s^h\Vert_{0,1}^2 + \mu\,\left[ \vert\f_s^h\vert_{1,-1}^2 + 2\,\Vert\f_s^h\Vert_{0,-3}^2 \right] \right)\,c_k^h+
\mu\,\left(\tilde u^k_h \mid p_s^h \right)= \left( f^k \mid p_s^h \right)\,. 
\label{axi.cmu5}
\end{eqnarray}
However, like any function in~$V^h_\circ$, $\f_s^h$ does not necessarily
belong to~$\HOneMOne$ or $L^2_{-3}(\om)$. This is no serious problem: 
like in the
prismatic case \cite[\S5.3]{CJK+04a}, we shall rather discretise the 
representation formula~(\ref{axi.oper3}), which we rewrite as follows:
\begin{equation}
c_k =\frac1{\Vert p_s\Vert_{0,1}^2}\,\left[ \left( f^k \mid p_s \right) - 
\mu\, \left( z^k \mid q_s \right) \right],
\label{axi.oper3b}
\end{equation}
where $q_s = p_s/r^2$ and $z_k=A_k^{-1}f^k\in\VOneOneZ$ solves  
\begin{equation}\label{axi.zmu}
a_k\left( z^k,v \right) = \left( f^k \mid v \right),\quad \forall 
v\in\VOneOneZ.
\end{equation} 
So, we state the 

\smallskip

\no\underline{{\bf SCM Algorithm} for finding 
$\tilde u^k_h\in V^h_\circ$ and $c_k^h\in \R^1$}.

\smallskip 

\no\underline{\bf Step 1}. Find $z^k_h\in V^h_\circ$ such that 
\begin{equation}\label{axi.zmuh}
a_k(z^k_h, v) = \left( f^k \mid v \right) \quad \forall\, v\in V^h_\circ\,.
\end{equation} 
Compute $c_k^h$ as follows:
\begin{equation}
c_k^h =\frac1{\Vert p_s^h\Vert_{0,1}^2}\,\left[ \left( f^k \mid p_s^h \right) - \mu\, \left( z^k_h \mid q_s^h \right) \right],
\quad \mbox{if} \quad k< C^\star\,h^{-\frac{1}{2-\alpha_0}}\,;
\label{axi.oper3a}
\end{equation}
for some fixed constant~$C^\star$, and 
\begin{equation}
c_k^h =0 \quad \mbox{if} \quad k\ge C^\star\,h^{-\frac{1}{2-\alpha_0}}\,.
\label{axi.oper3d}
\end{equation}

\no\underline{\bf Step 2}. Find  $\tilde u^k_h\in V^h_\circ$ such that 
\begin{equation} \label{axi.cmu6a}
a_k\left(\tilde u^k_h, v \right)
+c_k^h\,a_k\left(\f_s^h, v\right) 
 = \left( f^k \mid v \right) \quad \forall v\in V^h_\circ\,.
\end{equation}

\begin{lemma}\label{axi.lem:bound4a} For  the solution $z^k$ to 
the problem (\ref{axi.zmu}) and its piecewise linear finite element 
approximation $z^k_h$ in (\ref{axi.zmuh}), 
we have the following error estimates 
\begin{eqnarray}
\left\|z^k-z^k_h\right\|_{0,-1} &\ls& k^{-2}\,\|f^k\|_{0,1}\,, 
\label{axi.zmu5}\\ 
%\left\|z^k-z^k_h\right\|_{0,-1} &\ls& h^{\a_0}\,\|f^k\|_{0,1}\,, \nonumber\\
\left\|z^k-z^k_h\right\|_{0,-1} &\ls& k^{-1}\, \left[ h^{\alpha_0}\,k^{\alpha_0-1} + h \right]\,\|f^k\|_{0,1}\,,
\label{axi.zmu6}\\ 
\left\|z^k-z^k_h\right\|_{0,1} &\ls& \left[ h^{2\alpha_0}\,k^{2(\alpha_0-1)} + h^2 \right]\,
\|f^k\|_{0,1}\,, \label{axi.zmu4}
\end{eqnarray}
while for the coefficient $c_k$ in (\ref{axi.oper3b}) and its 
approximation $c_k^h$ in (\ref{axi.oper3a}), we have 
\begin{equation}\label{axi.error1} 
\left| c_k-c_k^h \right|
 \ls (h^{2\alpha_0}\,k^{2\alpha_0} + h^2\,k^2)\,\|f^k\|_{0,1} \,.
\end{equation}
\end{lemma}
\noindent\textit{Proof:} It follows from (\ref{axi.zmu}) and (\ref{axi.zmuh}) 
that 
\[
 a_k\left( z^k-z^k_h, z^k-z^k_h \right)=a_k\left( z^k, z^k-z^k_h \right)  = \left(f^k \mid z^k-z^k_h \right).
\]
This implies 
\[
 \left| z^k-z^k_h \right|^2_{1,1} + k^2\, \left\|z^k-z^k_h\right\|^2_{0,-1} 
  \le  \left\|f^k\right\|_{0,1}\,\left\|z^k-z^k_h\right\|_{0,1} \le \rmax\, \left\|f^k\right\|_{0,1}\,\left\|z^k-z^k_h\right\|_{0,-1},
\]
hence (\ref{axi.zmu5}).
Then, using C\'ea's lemma, Lemma~\ref{axi.lem:interpo:Vs}, Thm~7.1 of~\cite{MeRa82}, and the bounds (\ref{axi.eq:uk2}) and~(\ref{axi.asympu}), we obtain another estimate:
\begin{eqnarray*}
\left\| z^k - z^k_h \right\|^2_{(k)} &\le& \left\| z^k - P_h z^k \right\|^2_{(k)}\\
&=& \bigIII z^k - P_h z^k \bigIII_{1,1}^2 + (k^2-1)\, \left\| z^k - P_h z^k \right\|_{0,-1}^2\\
&\ls& h^{2\alpha_0}\, \left| z^k \right|_{1+\a_0,1}^2 + (k^2-1)\,h^2\, \left| z^k \right|_{1,-1}^2,\\
&\ls& \left[ h^{2\alpha_0}\,k^{2(\alpha_0-1)}  + h^{2} \right]\, \left\|f^k\right\|_{0,1}^2.
\end{eqnarray*}
Of course, a similar bound holds for any $g\in L{}^2_1(\om)$, $w=A_k^{-1}g$ and~$P_h w$. Thus, the estimate~(\ref{axi.zmu4}) follows from a duality argument like in~\cite[Lemma~5.3]{CJK+04a}.
Moreover, we obtain~(\ref{axi.zmu6}) thanks to the bound: $\|\cdot\|_{0,-1}^2 \le k^{-2}\, \|\cdot\|_{(k)}^2$.

\medbreak

To obtain the estimate~(\ref{axi.error1}), we subtract (\ref{axi.oper3b}) from~(\ref{axi.oper3a}) to obtain
\begin{eqnarray*}
c_k -c_k^h
&=& \Big\{ \frac{\left(f^k\mid p_s \right)}{\Vert p_s\Vert_{0,1}^2} - 
   \frac{\left(f^k\mid p_s^h \right)}{\Vert p_s^h\Vert_{0,1}^2}\Big\} +\mu\, 
  \Big\{ \frac{\left(z^k_h \mid q_s^h \right)}{\Vert p_s^h\Vert_{0,1}^2}-
 \frac{\left(z^k \mid q_s \right)}{\Vert p_s\Vert_{0,1}^2}\Big\}
\bydef I_1 + I_2.
\end{eqnarray*}
We bound $I_1$ by Lemma~\ref{axi.lem:bound2}: 
$|I_1|\ls h^{2\a_0}\,\left\|f^k\right\|_{0,1}$.
As for $I_2$, it is zero when $\mu=0$; otherwise we rewrite it as 
follows:
\begin{eqnarray*}
 \frac{I_2}{\mu} &=& 
\frac{1}{\left\| p_s^h \right\|_{0,1}^2 }\, \left\{ 
\left(  z^k_h - z^k \mid q_s \right) + 
\left( z^k_h - z^k \mid q_s^h - q_s \right) +
\left( z^k \mid q_s^h - q_s \right) \right\} \\
&&\mbox{} \hfill \mbox{}+ \left( z^k \mid q_s \right)\Big\{\frac{1}{\Vert p_s^h\Vert_{0,1}^2}  - \frac{1}{\Vert p_s\Vert_{0,1}^2}\Big\} \\
&\bydef& J_2^1 + J_2^2 + J_2^3 + J_2^4.
\end{eqnarray*}
Then, recalling that $\|q_s\|$ is constant, we estimate:
\begin{itemize}
\item $|J_2^1| \ls \left \| z^k_h - z^k \right\|_{0,1} \ls \left[ h^{2\alpha_0}\,k^{2(\alpha_0-1)}  + h^{2} \right]\, \left\|f^k\right\|_{0,1}$ by~(\ref{axi.zmu4}).
\item $|J_2^2| \ls h^{\a_0}\,\left\|z^k_h - z^k\right\|_{0,-1} \ls \left[ h^{2\alpha_0}\,k^{\alpha_0-2} + h^{1+\alpha_0}\,k^{-1}  \right]\,\|f^k\|_{0,1}$ by (\ref{axi.approx:tildeq:0.-1.w}) and~(\ref{axi.zmu6}).
\item $|J_2^3| \ls h^{2\a_0}\,\left\|z^k\right\|_{0,-3} \ls h^{2\a_0}\,k^{-2}\,\left\|f^k\right\|_{0,1}$ by (\ref{axi.approx:tildeq:0.-3.w}) and~(\ref{axi.eq:uk2}).
\item
$|J_2^4| \ls h^{2\a_0}\,\left\|z^k\right\|_{0,1} \le h^{2\a_0}\,\rmax\, \left\|z^k\right\|_{0,-1} \ls h^{2\a_0}\,k^{-2}\,\left\|f^k\right\|_{0,1}$ by (\ref{axi.eq:norme:pse:0.1}) and~(\ref{axi.eq:uk1}).
\end{itemize}
Summarising, we obtain
\begin{eqnarray*}
\left| c_k - c_k^h \right| &\le& \left| I_1 \right| + \mu\,\left\{  \left| J_2^1 \right| + \left| J_2^2 \right| + \left| J_2^3 \right| + \left| J_2^4 \right| \right\} \\ 
&\ls& \left( h^{2\a_0}  + h^{2\a_0}\,k^{\a_0} + h^{2\a_0}\,k^{2\a_0} + h^{1+\a_0}\,k + h^2\,k^2 \right)\, \left\|f^k\right\|_{0,1}.
\end{eqnarray*}
The estimate~(\ref{axi.error1}) then follows by remarking that the first,
second and fourth terms in the bracket are negligible with respect to 
the third.\hfill$\Box$

\bigbreak

Now, we observe that the formula~(\ref{axi.cmu6a}) for computing~$\tilde u^k_h$, as well as the SCM reconstruction formula for $u^k_h$
\begin{equation}\label{axi.uh1}
u^k_h=\tilde u^k_h +c_k^h\,\f_s^h 
 = \tilde u^k_h +c_k^h\,(\tilde \f_h+\delta_h \f_{_P}).  
\end{equation}
are formally identical to their prismatic counterparts 
(cf.~the SCM algorithm of~\cite[\S5.3]{CJK+04a}); 
and the ``building blocks'' $c_k^h$ and~$\f_s^h$ also 
satisfy estimates similar to those of the prismatic case. Indeed, under the
assumption $k< C^\star\,h^{-\frac{1}{2-\alpha_0}}$, both terms within the
bracket in~(\ref{axi.error1}) are negligible with respect to~$h\,k$.
Hence the following two results, whose proofs closely parallel that of 
Lemma~5.4 and Theorem~5.1 in~\cite{CJK+04a}, with the same kind of 
adaptations (use of weighted norms, $P_h$ and 
Lemma~\ref{axi.lem:interpo:Vs}) as usual.

\begin{lemma}\label{axi.lem:bound4} The following error estimate holds:
\[ %\begin{equation}\label{axi.umu2a}
\left \|\tilde u^k - \tilde u^k_h \right\|_{(k)}^2 \ls k\,
\left( h^2\,(1+k^2\,h^2)\, \left\|f^k\right\|_{0,1}^2 + \left| c_k-c_k^h 
\right|^2 \right).
\] %\end{equation}
\end{lemma}

\begin{theorem}\label{axi.thm:main1}
% Assume that $N\,h\in[\lambda_1,\lambda_2]$, with $0<\lambda_1<\lambda_2$.
Let $u^k$ be the solution to the equation (\ref{axi.laplace1:Delta.k}--\ref{axi.laplace1:cl})
and $u^k_h$ be its finite element approximation given in (\ref{axi.uh1}). 
Then the following error estimate holds:
\begin{equation}
\left\| u^k-u^k_h \right\|_{(k)} \ls k^2\,h\,\left\|f^k\right\|_{0,1}
\label{axi.err:k>2}
\end{equation}
\end{theorem}

\subsection{Approximation of the singular functions for the modes 
$|k|=0,\ 1$.\label{axi.sec:sg:petits:modes}}
The FE approximation of these functions has been exposed 
in~\cite[\S\S4.1 and~4.2]{ACLS03}. (In that work, the Laplacians $\Delta_0$ 
and~$\Delta_1$ are respectively called $\Delta$ and~$\Delta'$). 
We keep this method, with the following modification.
The dual singular functions associated to the reentrant edge undergo the
same inconvenience as $p_s^2$, namely, the Laplacian of the principal
parts as defined in~\cite{ACLS03} do not belong to~$L{}^2_1(\om)$. 
Hence we must enrich them, just as we did for $p_{_P}^2$, in order to
preserve the convergence rate. Calculating like in~\S\ref{axi.sec:ps},
we obtain the following decompositions:
\begin{eqnarray}
p_s^1 = p_p^1 + \widehat p^1, && p_p^1 \bydef 
\rho^{-\a}\,\sin(\a\phi)\, \frac{r}a\, \left[ 1 - 
\frac{3\rho}{2a}\,\cos\phi' \right],
\quad \widehat p^1\in \VOneOne \,;\\
p_s^{0,e} = p_p^{0,e} + \widehat p^{0,e}, && p_p^{0,e} \bydef 
\rho^{-\a}\,\sin(\a\phi)\, \left[ 1 - \frac{\rho}{2a}\,\cos\phi' \right],
\quad \widehat p^{0,e}\in \HOneOne.
\end{eqnarray}
The Laplacians of the principal parts are:
\begin{eqnarray}
\vt^1_p \bydef \Delta_1 p^1_p &=& \frac{1}{a^2}\,\rho^{-\a}\, 
\left[ -\frac{9}{2}\,\sin(\a\phi) + \frac{15}{2}\,\a\,\cos\phi'\,
\sin(\a\phi + \phi')  \right] \in L{}^2_1(\om) \,;\qquad\\
\vt^{0,e}_p \bydef \Delta_{0} p^{0,e}_p &=& \frac{1}{a\,r}\,\rho^{-\a}\, \left[ -\frac{1}{2}\,\sin(\a\phi) + \frac{3}{2}\,\a\,\cos\phi'\,\sin(\a\phi + \phi')  \right] \in L{}^2_1(\om).\qquad
\end{eqnarray}
Then we proceed like in~\S\ref{axi.sec:ps} to obtain:
\begin{eqnarray}
%\mbox{Mode 1:}\qquad
\BigIII p_s^1-p_s^{1;h} \BigIII_{1,1} \ls h^{\alpha_0}\,,&&
\left\| p_s^1-p_s^{1;h} \right\|_{0,1} \ls h^{2\alpha_0}\,, \\
%\mbox{Mode 0:}\qquad
\left| p_s^{0,e}-p_s^{0,e;h} \right|_{1,1} \ls h^{\alpha_1}\,,&&
\left\| p_s^{0,e}-p_s^{0,e;h} \right\|_{0,1} \ls h^{2\alpha_1}\,.
\end{eqnarray}
However,
for the primal edge singular functions, the method of~\cite{ACLS03} 
yields the desired convergence rate. We just recall the decompositions:
\begin{equation}
\f_s^{k} = \tilde\f^{k} + \delta^{k}\,\f_{_P}^{k},\quad \tilde\f^{k} \in 
H{}^2_1(\om)\cap H^1_{(k)}(\om),\quad \f_{_P}^{k} = \left( \frac{r}{a} 
\right)^k\, \rho^{\a}\,\sin(\a\phi),
%\label{axi.gamma1}
\end{equation}
as well as the Laplacians of the principal parts:
\begin{equation}
\psi_{_P}^{k}\bydef\Delta_k\f_{_P}^{k} = 
\frac{(k+1)(k+2)}{2\,a^k\,r^{1-k}}\, \a\,\rho^{\a-1}\,
\sin\left[(\a-1)\,\phi - \phi_0 \right].
\end{equation}
The line of proof already exposed in~\S\ref{axi.sec:phis} then easily 
leads to the error estimates:
\begin{equation}
\bigIII \f_s^1 - \f_s^{1;h} \bigIII_{1,1} \ls h \quad\mbox{and}\quad
\bigl|  \f_s^{0,e} - \f_s^{0,e;h} \bigr|_{1,1} \ls h.
\end{equation}

\medbreak

Now, as far as the conical point singularities are concerned, the 
method appears very similar to that of \cite[\S\S5.1 and~5.2]{CJK+04a} 
since the principal parts $p_{_P}^{0,c}$ and $\f_{_P}^{0,c}$
have a vanishing Laplacian~$\Delta_0$.
So, \emph{mutatis mutandis}, we get the error estimates:
\begin{equation}
\left| p_s^{0,c} - p_s^{0,c;h} \right|_{1,1} \ls h^{\a_1},\quad  
\left\| p_s^{0,c} - p_s^{0,c;h} \right\|_{0,1} \ls h^{2\a_1}, \quad
\left| \f_s^{0,c} - \f_s^{0,c;h} \right|_{1,1} \ls h.
\end{equation}

\smallbreak

\begin{remark}
Thanks to the asympotic expansions~\cite[Eq.~8.7.1]{AbSt65} of the
Legendre function, it is possible to compute the function~$P_\nu(\cos\phi)$ 
with an arbitrary precision. Thus, one can compute \emph{once and for all}
the singularity exponent~$\nu$ and the integral 
in~(\ref{eq:formule:delta.0.c}) with an accuracy equal to the machine 
precision. All this guarantees that the errors due to the approximation of 
the conical singular functions will be  negligible before the FE 
discretisation error. 
\end{remark}

\subsection{Approximation of $\tilde u^k$ and $c_k$, for $|k|\le1$.} 
As the representation formulae (\ref{axi.oper3:k=1}) and~(\ref{axi.oper3:k=0})
for the singularity coefficients of these modes are rather standard, one
can use the simple discrete versions:
\begin{equation}
c_{\pm1}^h =\frac{\left( f^{\pm1}  \vert \,p_s^{1;h} \right) }{\Vert p_s^{1;h}\Vert_{0,1}^2}\,,\qquad
c_{0,j}^h =\frac{\left( f^0 \mid p_s^{0,j;h} \right)}{\Vert p_s^{0,j;h}\Vert_{0,1}^2}\,.
%\label{axi.oper3:k=1}
\end{equation}
Similarly, we will approximate the regular parts $\tilde u^k,\ |k|\le1$, 
by $\tilde u^1_h,\ \tilde u^{-1}_h\in V^h_\circ$, and $\tilde u^0_h \in 
V^h_\diamond$ such that
\begin{eqnarray}
k=\pm1 :&& \negthickspace\negthickspace  
a_1\left(\tilde u^k_h, v_h \right) +c_k^h\,a_1\left(\f_s^{1,h}, v_h\right)
=\left( f^k \mid v_h \right), \quad \forall v_h\in V^h_\circ\,,\\
k=0 :&& \negthickspace\negthickspace 
a_0\left( \tilde u^0_h,v_h \right) + c_{0,e}^h\,a_0\left( \f_s^{0,e;h}, v_h\right) + c_{0,c}^h\,a_0\left( \f_s^{0,c;h}, v_h \right) = \left( f^0 \mid v_h\right),\quad\forall v_h\in V^h_\diamond\,.
\qquad\quad
\end{eqnarray}
Of course, we have the SCM reconstruction formulae:
\begin{eqnarray}
k=\pm1 :&& u^k_h=\tilde u^k_h +c_k^h\,\f_s^{1;h} 
 = \tilde u^k_h +c_k^h\,(\tilde \f^1_h+\delta^1_h \f^1_{_P}) \,;
\label{axi.uh1:k=1}\\
k=0:&& u^0_h=\tilde u^0_h + c_{0,e}^h\,\f_s^{0,e;h} + c_{0,c}^h\,\f_s^{0,c;h}
\nonumber\\
&& \hphantom{u^0_h}= \tilde u^0_h + c_{0,e}^h\,(\tilde \f^{0,e}_h+\delta^{0,e}_h \f^{0,e}_{_P})
+ c_{0,c}^h\,(\tilde \f^{0,c}_h+\delta^{0,c}_h \f^{0,c}_{_P}).
\label{axi.uh1:k=0}
\end{eqnarray}
The results of \S\ref{axi.sec:sg:petits:modes} then allow to conclude that:
\begin{eqnarray}
k=\pm1:&&
\left| c_k - c_k^h \right| \ls h\, \left\|f^k\right\|_{0,1},\quad 
\bigIII \tilde u^k - \tilde u^k_h \bigIII_{1,1} \ls h\,
\left\|f^k\right\|_{0,1}, \nonumber\\
&& \bigIII u^k - u^k_h \bigIII_{1,1} \ls h\,\left\|f^k\right\|_{0,1}  \,;
\label{axi.err:k=1}
\\
k=0:&&
\left| c_{0,j}^h - c_{0,j} \right| \ls h\,\left\|f^0\right\|_{0,1},
\quad \big| \tilde u^0 - 
\tilde u^0_h \big|_{1,1} \ls h\,\left\|f^0\right\|_{0,1}, \nonumber\\
&&\big| u^0 - u^0_h \big|_{1,1} \ls h\,\left\|f^0\right\|_{0,1}.
\label{axi.err:k=0}
\end{eqnarray}

\section{Fourier Singular Complement Method\label{axi.sec:FSCM}}
Let $u$ be the solution to the 3D problem~(\ref{axi.Poisson}), and
$u^k$ its Fourier coefficients. From the previous Sections, 
we know that $u^k(r,z)$ solves the 2D problem~(\ref{axi.Poissonk}), 
the weak formulation of the elliptic 
problem~(\ref{axi.laplace1:Delta.k}--\ref{axi.laplace1:cl}).
And, according to the mode~$k$, one can decompose $u^k$ 
as~(\ref{axi.observ1}), (\ref{axi.observ1:k=1}) 
or~(\ref{axi.observ1:k=0}).

\medbreak

The result of Heinrich~{\cite[Thm~5.2]{Hein93}} can be straightforwardly
extended to our domain with a sharp vertex.
\begin{theorem}
\label{axi.lem:f6}
Let $f\in h^2(\Om)$, and $u\in \overset{\circ}{H}{}^1(\Om)$ be the solution to~(\ref{axi.Poisson}). Then:
\begin{equation}\label{axi.splitting2}
u(r,\t,z)=\tilde u(r,\t,z) + \ga(\t)\,\f_s^{2}(r,z) + 
c_0^c\,\,\f_s^{0,c}(r,z),
\end{equation}
with: $\tilde u\in H^2(\Om)\cap \overset{\circ}{H}{}^1(\Om)$, and
$\ga\in H^2(\mathbf{S}^1)$ is given by the formula:
$$\ga(\t) = \frac{\delta^{0,e}}{\delta^2}\,c_0^e +
\frac{\delta^{1}}{\delta^2}\,\sum_{k=\pm1} c_k\, \E^{\I k \theta}
+ \sum_{|k|\ge2} c_k\, \E^{\I k \theta}.$$
\end{theorem}
Like in the prismatic case (cf.~\cite{CJK+04a}, Remark~6.1), the hypothesis
$f\in h^2(\Om)$ is crucial: the lack of its satisfaction would prevent the convergence of~$\ga$ in a regular enough space, and hence
that of the singular part of the solution in the natural space.

\medbreak

We define the Fourier--SCM (FSCM) solution to (\ref{axi.Poisson}) as follows: 
\[ %\begin{equation} \label{axi.scm2}
 u_{h}^{[N]}  =\sum_{k=-N}^N u^k_h(r,z) \,\E^{\I k\theta},
\] %\end{equation}
where $u^k_h$ is the SCM solution to~(\ref{axi.Poissonk}) algorithmically
defined in \S\ref{axi.sec:SCM}.
The main result on this method is the following
\begin{theorem}
Assume $f\in h^2(\Om)$. Then
the following error estimate holds: 
\[
\left|  u- u_{h}^{[N]}\right|_{H^1(\Om)} \ls 
(h+N^{-1})\, \Big\{ \Big\|f\Big\|_{L^2(\Om)}+
%\Big\|\frac{\pa f}{\pa\t}\Big\|_{L^2(\Om)} + 
\Big\|\frac{\partial^2 f}{\partial\t^2}\Big\|_{L^2(\Om)} \Big\}\,.
\]
\end{theorem}
\noindent {\it Proof:} 
Using the definition of~$u_{h}^{[N]}$ and~(\ref{axi.eq:fourier:H1})
we have
\begin{displaymath}
\left|  u- u_{h}^{[N]}\right|_{H^1(\Om)}^2 =
\sum_{|k|\le N} \left\| u^k - u^k_h \right\|_{(k)}^2 + \sum_{|k|>N} 
\left\| u^k \right\|_{(k)}^2 \bydef E_1 + E_2.
\end{displaymath}
Using~(\ref{axi.eq:uk1}), we estimate~$E_2$ as:
\begin{displaymath}
E_2 \le N^{-2}\, \sum_{|k|>N} k^2\, \left( \left|u^k\right|_{1,1}^2 +
k^2\, \left\|u^k\right\|_{0,-1}^2 \right) \ls N^{-2} \sum_{|k|>N}
\left\| f^k \right\|_{0,1}^2 \le N^{-2}\, \|f\|_{L^2(\Om)}^2 \,.
\end{displaymath}
As for~$E_1$, we cut it into three parts, corresponding to $k=0$, $|k|=1$, 
and $2\le|k|\le N$, which we bound respectively by (\ref{axi.err:k=0}) 
and (\ref{axi.err:k=1}) and~(\ref{axi.err:k>2}):
\begin{eqnarray*}
E_1 &\ls& h^2\, \left\|f^0\right\|_{0,1}^2 + h^2\, \left( 
\left\|f^1\right\|_{0,1}^2 + \left\|f^{-1}\right\|_{0,1}^2 \right) + h^2\,
\sum_{2\le|k|\le N} k^4\, \left\|f^k\right\|_{0,1}^2\\
&\ls& h^2\, \left\{ \|f\|_{L^2(\Om)} + 
\left\|\frac{\partial^2 f}{\partial\t^2}\right\|_{L^2(\Om)} \right\}\,,
\end{eqnarray*}
where we have used Lemma~\ref{axi.lem:completeness} to bound the sum. 
Hence the result.\hfill$\Box$

\section{Conclusion}
In this paper, we have proven that the FSCM for the Poisson equation
achieves the optimal convergence rate for $\mathbb{P}_1$ finite elements
and a datum of $L^2$-style regularity in the meridian directions.
The same result also holds for the discretization of the Poisson problem with a homogeneous Neumann boundary condition, or with non-homogeneous boundary conditions, provided there exist sufficiently smooth liftings. 

\medbreak

This result closely parallels that of the companion paper~\cite{CJK+04a}.
The specificities of the axisymmetric geometry (namely, that the 
2D~problems are set in weighted Sobolev spaces, which moreover vary for the
low-order Fourier modes before stabilising, and involve differential 
operators with non-constant coefficients) only cause technical difficulties.
As far as the presence of conical vertices is concerned, its effect is 
no more than a finite-dimensional perturbation. 
Furthermore, it is no difficulty to consider the case of an axisymmetric 
domain $\Om$ with several reentrant edges (i.e. $\om$ with  several 
off-axis reentrant corners) and/or several sharp vertices.

\medbreak

As already mentioned, this paper is the second part of a three-part article \cite{CJK+04a,CJK+04c}.
In~\cite{CJK+04c}, the FSCM is analysed from a numerical point of view (complexity, implementation issues, numerical experiments, etc.), and it is compared to other methods---in the axisymmetric case, to anisotropic mesh refinement techniques.

\medbreak

One can apply the same theoretical and numerical techniques to the fully
axisymmetric heat or wave equations, with any $L^2$-smooth (in space) right-hand side. For these PDEs, the 
singular functions $p_s$ and $\f_s$ do not depend on the time-step. \\
Finally, the results, can also be viewed as the first effort towards the discretization of electromagnetic 
fields in axisymmetric domains, with {\em continuous} numerical approximations, the importance of which is well-known,
cf.~\cite{ADHR93}. As a matter of fact, the SCM developed in~\cite{AsCL02,AsCL03,ACLS03} for fully axisymmetric
electromagnetic computations can be generalized to arbitrary data, with the help of the results obtained here.

\end{document}